\documentclass[]{interact}

\usepackage{epstopdf}
\usepackage[caption=false]{subfig}

\usepackage{url,hyperref}
\usepackage{algorithm}
\usepackage[usenames,dvipsnames]{xcolor}
\usepackage{algpseudocode}
\usepackage{algorithmicx}
\usepackage[ruled,vlined, linesnumbered, algo2e]{algorithm2e}
\usepackage[numbers,sort&compress]{natbib}
\bibpunct[, ]{[}{]}{,}{n}{,}{,}

\theoremstyle{plain}
\newtheorem{theorem}{Theorem}[section]
\newtheorem{lemma}[theorem]{Lemma}

\theoremstyle{definition}

\theoremstyle{remark}
\newtheorem{remark}{Remark}

\newtheorem{assumption}{Assumption}[section]

\newcommand{\Id}{\mathbb{I}}
\newcommand{\R}{\mathbb{R}}
\newcommand{\Rext}{\R\cup\{+\infty\}}
\newcommand{\set}[1]{\left\{#1\right\}}
\newcommand{\sets}[1]{\{#1\}}
\newcommand{\norm}[1]{\left\Vert#1\right\Vert}
\newcommand{\norms}[1]{\Vert#1\Vert}
\newcommand{\Eproof}{\hfill $\square$}
\newcommand{\prox}{\mathrm{prox}}

\newcommand{\relint}[1]{\mathrm{ri}\left(#1\right)}
\newcommand{\argmin}{\mathrm{arg}\!\displaystyle\min}
\newcommand{\dom}[1]{\mathrm{dom}(#1)}
\newcommand{\Exps}[2]{\mathbb{E}_{#1}\left[#2\right]}
\newcommand{\Prob}[1]{\mathbf{Prob}\left(#1\right)}
\newcommand{\iprod}[1]{\left\langle #1\right\rangle}
\newcommand{\iprods}[1]{\langle #1\rangle}
\newcommand{\Exp}[1]{\mathbb{E}\left[#1\right]}

\newcommand{\BigO}[1]{\mathcal{O}\left(#1\right)}

\newcommand{\BigOs}[1]{\mathcal{O}(#1)}
\newcommand{\BigOc}[1]{\mathcal{O}\left(#1\right)}
\newcommand{\Up}[1]{\mathbb{U}_{\mathbf{q}}\left(#1\right)}

\newcommand{\Xc}{\mathcal{X}}
\newcommand{\Yc}{\mathcal{Y}}
\newcommand{\Zc}{\mathcal{Z}}
\newcommand{\Wc}{\mathcal{W}}
\newcommand{\Lc}{\mathcal{L}}
\newcommand{\Tc}{\mathcal{T}}
\newcommand{\Sc}{\mathcal{S}}
\newcommand{\Fc}{\mathcal{F}}
\newcommand{\Gc}{\mathcal{G}}

\newcommand{\Cc}{\mathcal{C}}
\newcommand{\Ec}{\mathcal{E}}

\newcommand{\mytbi}[1]{\textbf{\textit{#1}}}
\newcommand{\mytb}[1]{\textbf{#1}}

\newcommand{\beforesec}{\vspace{-2ex}}
\newcommand{\aftersec}{\vspace{-2ex}}
\newcommand{\beforesubsec}{\vspace{-3ex}}
\newcommand{\aftersubsec}{\vspace{-1ex}}
\newcommand{\beforesubsubsec}{\vspace{-3ex}}
\newcommand{\aftersubsubsec}{\vspace{-1ex}}
\newcommand{\beforepara}{\vspace{-3ex}}

\newcommand{\myeq}[2]{ 
\vspace{-0.5ex}
\begin{equation}\label{#1}
{#2}
\vspace{-0.5ex}
\end{equation}
}
\newcommand{\myeqn}[1]{ 
\vspace{-0.5ex}
\begin{equation*}
{#1}
\vspace{-0.5ex}
\end{equation*}
}

\newcommand{\rvtext}[1]{#1}

\begin{document}


\title{A New Randomized Primal-Dual Algorithm for Convex Optimization with Optimal Last Iterate Rates}

\author{
\name{Quoc Tran-Dinh\textsuperscript{a}\thanks{CONTACT: Quoc Tran-Dinh. Email: quoctd@email.unc.edu} and Deyi Liu\textsuperscript{a}}
\affil{\textsuperscript{a}Department of Statistics and Operations Research\\ 
The University of North Carolina at Chapel Hill\\ 
318 Hanes Hall, UNC-Chapel Hill, NC 27599}
}

\maketitle

\begin{abstract}
We develop a novel unified randomized block-coordinate primal-dual algorithm to solve a class of nonsmooth constrained convex optimization problems, which covers different existing variants and model settings from the literature.
We prove that our algorithm achieves optimal  $\BigO{n/k}$ and $\BigO{n^2/k^2}$ convergence rates (up to a constant factor) in two cases: general convexity and strong convexity, respectively, where $k$ is the iteration counter and $n$ is the number of block-coordinates.
Our convergence rates are obtained through three criteria: primal objective residual and primal feasibility violation, dual objective residual, and primal-dual expected gap.
Moreover, our rates for the primal problem are on the last iterate sequence.
\rvtext{Our dual convergence guarantee requires additionally a Lipschitz continuity assumption.}
We specify our algorithm to handle two important special cases, where our rates are still applied. 
Finally, we verify our algorithm on two well-studied numerical examples and compare it with two existing methods.
Our results show that the proposed method has encouraging performance on different experiments.
\end{abstract}


\beforesec
\section{Introduction}\label{sec:intro}
\aftersec
We consider the following nonsmooth constrained convex optimization problem:
\myeq{eq:primal_form}{
F^{\star} := \min_{x \in\R^p, w \in \R^\rvtext{m}}\Big\{ F(x, w) :=  h(x) + \sum_{i=1}^nf_i(x_i) + g(w) \quad \text{s.t.} \quad Kx + Bw = b \Big\},  \tag{P}
}
where $f_i : \R^{p_i} \to\Rext$ for $i=1,\cdots, n$, are proper, closed, possibly nonsmooth, and convex functions, $p_1 + \cdots + p_n = p$, $h : \R^p\to\R$ is a smooth and convex function, 
$g : \R^m \to \Rext$ is a proper, closed, possibly nonsmooth, and convex function, $K\in\R^{d\times p}$, $B \in \R^{d\times m}$, and $b\in\R^d$ are given. 
For notational simplicity, let us denote by $f(x) := \sum_{i=1}^nf_i(x_i)$ throughout this paper.

Note that \eqref{eq:primal_form} looks simple, but it is sufficiently general to cope with a broad class of convex optimization problems in practice, ranging from unconstrained to constrained settings, including the composite models considered in \cite{alacaoglu2019convergence,Alacaoglu2017,chambolle2017stochastic,tan2018stochastic,zhang2017stochastic}.
In particular, \eqref{eq:primal_form} also covers conic programming (e.g., linear, convex quadratic, second-order cone, and semidefinite programming), image and signal processing, machine learning, network and distributed optimization, and optimal transport, see, e.g., \citep{Boyd2004,Bertsekas1989b,rockafellar1984network,TranDinh2012h}.

The corresponding dual problem of \eqref{eq:primal_form} can be written as
\myeq{eq:dual_form}{
D^{\star} := \max_{y \in\R^d} \ D(y), 
~~\text{where}~
D(y) := \min_{x, w}\big\{ F(x, w) + \iprods{Kx + Bw - b, y} \big\}.
\tag{D}
}
Here, $D$ is called the dual function of \eqref{eq:primal_form}.
Our goal in this paper is to develop a new unified randomized block-coordinate primal-dual algorithm to simultaneously solve both \eqref{eq:primal_form} and its dual \eqref{eq:dual_form} which is simple to implement and can achieve state-of-the-art convergence rates without imposing strong assumptions on \eqref{eq:primal_form}.

\vspace{1ex}
\noindent\textbf{Motivation.}
We are interested in the case when the primal dimension $p$ of $x$ is sufficiently large so that computing the full gradient of $h$ and the proximal operator of $f$ can be prohibited. 
However, the dimensions $m$ of $w$ and $d$ of $y$ are relatively small, so that full operations on these spaces can be computed efficiently.
This structure is sufficiently generic to cope with many existing models, including \cite{alacaoglu2019convergence,Alacaoglu2017,chambolle2017stochastic,tan2018stochastic,zhang2017stochastic}. 

Although randomized primal-dual algorithms for solving \eqref{eq:primal_form} have been widely studied in the literature, including \citep{konevcny2017semi,tan2018stochastic,yu2015doubly}, it remains unclear if one can achieve an ``optimal'' convergence rate on the last iterate under only convexity or strong convexity. 
Moreover, many works only focus on the unconstrained setting with smoothness assumption as opposed to the nonsmooth setting \eqref{eq:primal_form}.
For example, \citep{konevcny2017semi,yu2015doubly} only consider the unconstrained case when the objective function is both strongly convex and $L$-smooth to achieve a linear convergence rate, while \citep{tan2018stochastic} only achieves $\BigOs{1/\sqrt{k}}$ convergence rates under convexity, where $k$ is the iteration counter.
Our goal is to combine the augmented Lagrangian framework, randomized block-coordinate strategy \citep{fercoq2015accelerated}, and Nesterov's accelerated scheme \citep{Nesterov1983} to develop a new algorithm, Algorithm~\ref{alg:A1_main}, with optimal $\BigOs{n/k}$ convergence rate (up to a constant factor).
This convergence rate is much faster than $\BigOs{n/\sqrt{k}}$ in \citep{tan2018stochastic}.
Moreover, compared to other existing methods in \citep{Alacaoglu2017,chambolle2017stochastic,shalev2014accelerated}, Algorithm~\ref{alg:A1_main} achieves convergence rates on the last primal iterate $(x^k, w^k)$.
These rates can be boosted up to optimal $\BigO{n^2/k^2}$ rates (up to a constant factor) on the last primal iterate $(x^k, w^k)$ when either $f$ or $h$ is strongly convex.

\vspace{1ex}
\noindent\textbf{Related work.}
Solution methods for solving \eqref{eq:primal_form} and \eqref{eq:dual_form} have attracted great attention in recent years, including penalty schemes, augmented Lagrangian frameworks, primal-dual hybrid gradient method (PDHG), proximal splitting algorithms,  and variational inequality tools, see, e.g.,  \citep{Bauschke2011,BricenoArias2011,Chambolle2011,Chen2013a,combettes2012primal,combettes2011,Esser2010,Facchinei2003,Goldstein2013,He2012b,tran2017proximal}, just to name a few.
Notably, in \citep{Esser2010}, the authors generalize PDHG and show that their framework covers proximal forward-backward splitting (PFBS), the alternating direction
method of multipliers (ADMM), and the Douglas–Rachford splitting as special cases.
Hithereto, first-order primal-dual methods are perhaps the most popular ones for solving \eqref{eq:primal_form}, see, e.g., \citep{Bauschke2011,chambolle2016introduction,Esser2010a}.
In terms of convergence analysis, both linear convergence rate under strong convexity and smoothness assumptions and sublinear convergence rate under weaker assumptions are well established for primal-dual methods, including \citep{Bot2012,chambolle2016ergodic,Davis2014,Davis2014b,He2012b,Monteiro2010,TranDinh2015b}.  
However, most existing convergence rates are achieved via averaging or weighted averaging sequences instead of the last iterate sequence.

There also exist many stochastic primal-dual variants, including \citep{alacaoglu2020random,Alacaoglu2017,chambolle2017stochastic,fang2017faster,Lan2012,Ouyang2013,shalev2014accelerated,tan2018stochastic,zhang2017stochastic}, for solving \eqref{eq:primal_form} and its special instances.
For instance, one of the most notable works is \cite{chambolle2017stochastic}, which extends PDHG to stochastic variants.
As shown in  \cite{chambolle2017stochastic}, such stochastic variants often outperformance their deterministic counterparts.
As another example, \citep{zhang2017stochastic} also extends PDHG to a stochastic primal-dual coordinate variant, called SPDC.
This  method can achieve a linear convergence rate under the strong convexity and smoothness assumptions. 
In \citep{alacaoglu2019convergence,Alacaoglu2017}, the authors consider other stochastic variants of PDHG and obtain an optimal $\BigO{n/k}$-rate without strong convexity assumption.
Alternatively, stochastic alternating direction methods of multipliers (ADMM) have also been proposed to solve \eqref{eq:primal_form}, see, e.g., \cite{chambolle2017stochastic,fang2017faster,Ouyang2013,shalev2014accelerated,yu2017fast}.
Sublinear convergence rates in expectation or high probability have also been investigated for stochastic ADMM  \citep{chambolle2017stochastic,fang2017faster,shalev2014accelerated}. 
Recently, \citep{chambolle2017stochastic,fang2017faster,tan2018stochastic,zhang2017stochastic} show that stochastic primal-dual algorithms can perform several times faster than their deterministic counterparts for solving large-scale applications in machine learning. 
Other randomized block-coordinate methods and their asynchronous variants have been recently extended to monotone inclusions and general convex-concave minimax problems such as \cite{combettes2018asynchronous,combettes2015stochastic,peng2016arock,hamedani2018iteration}, which cover \eqref{eq:primal_form} as a special case. 
However, these algorithms are not accelerated and have no or slower non-ergodic convergence rates.
Recently, \cite{bai2020inexact,zhu2020stochastic} also extend other primal-dual methods to randomized and stochastic variants.
A very recent survey \cite{dvurechensky2021first} provides an excellent source on first-order methods for solving \eqref{eq:primal_form}, including randomized methods.

\vspace{1ex}
\noindent\textbf{Contribution.}
Our main contribution in this paper can be summarized as follows.
\begin{itemize}
\item[(a)]
We develop a unified randomized block-coordinate  primal-dual algorithm, Algorithm \ref{alg:A1_main}, to solve both \eqref{eq:primal_form} and \eqref{eq:dual_form}.
We prove optimal $\BigOc{n/k}$ convergence rates (up to a constant factor) for both \eqref{eq:primal_form} and \eqref{eq:dual_form} on three criteria: primal objective residual and primal feasibility violation, dual objective residual, and primal-dual expected gap, under only convexity and strong duality.
Moreover, our rates are on the last primal iterate $(x^k, w^k)$ compared to \cite{chambolle2017stochastic}.
\rvtext{Our dual convergence guarantee requires additionally the Lipschitz continuity of the conjugate $(f+h)^{\ast}$.}

\item[(b)]
If, in addition, $f$ or $h$ is strongly convex, then by appropriately adapting the parameter update rules, Algorithm \ref{alg:A1_main} can be boosted up to $\BigO{n^2/k^2}$ rates under the same three criteria. 
Again, our rates are optimal (up to a constant factor) and on the last primal iterate $(x^k, w^k)$.

\item[(c)]
We specify our algorithm to handle two special cases commonly studied in the literature: nonsmooth convex minimization with linear constraints and composite convex minimization.
In both cases, our convergence rates remain applicable.
\end{itemize}

\noindent\textbf{Comparison.}
Let us highlight the following aspects of our contribution.
Firstly, Algorithm~\ref{alg:A1_main} handles a more general class of problems than the composite model in, e.g., \cite{Alacaoglu2017,chambolle2017stochastic}.
Moreover, it is fundamentally different from SMART-CD in \cite{Alacaoglu2017}, where it  updates two dual variables and relies on the augmented Lagrangian framework instead of a smoothing technique as in \cite{Alacaoglu2017}. 
Algorithm~\ref{alg:A1_main} is also different from SPDHG in \cite{chambolle2017stochastic} since it is based on Nesterov's accelerated scheme with additional momentum steps.
Secondly,  if $g = 0$ and $B = 0$, then Algorithm~\ref{alg:A1_main} is also a fully randomized block-coordinate variant w.r.t. the primal variable $x$ as the one in \cite{tan2018stochastic}.
However, Algorithm~\ref{alg:A1_main} is accelerated.
It can be reduced to a non-accelerated variant as a special case of the method in \cite{tan2018stochastic}, see Subsection~\ref{subsubsec:special_case1}.
Thirdly,  by eliminating $w$, the dual update of Algorithm~\ref{alg:A1_main} possesses a three-point momentum step and uses dynamic parameter updates without any tuning.
This leads to a new type of algorithm, called ``non-stationary'' method \citep{liang2017local}.
Note that analyzing the convergence of ``non-stationary'' algorithms is often more challenging than that of stationary counterparts \citep{liang2017local}.
Fourthly, we establish three types of convergence guarantees, while most existing works only consider one.
\rvtext{Finally, compared to \cite{tran2019non}, both \cite{tran2019non} and Algorithm~\ref{alg:A1_main} exploit an augmented Lagrangian approach, and combining it with Nesterov's accelerated steps.
However, Algorithm~\ref{alg:A1_main} is different from \cite{tran2019non} on several aspects. 
First, \eqref{eq:primal_form} has different structure than the problem in  \cite{tran2019non} (see Subsection~\ref{subsubsec:special_case2}).
Second, Algorithm~\ref{alg:A1_main} is a randomized method, while  \cite{tran2019non} is deterministic. 
Third, it relies on the accelerated scheme in \cite{tseng2008accelerated}, while \cite{tran2019non} exploits Nesterov's original scheme in \cite{Nesterov1983}.
Fourth, Algorithm~\ref{alg:A1_main} unifies two cases (convex and strongly convex) in one single algorithm, while the second algorithm in \cite{tran2019non} for the strongly convex case is different and requires two proximal operations of $f$ per iteration.
Finally, the analysis of Algorithm~\ref{alg:A1_main} is much more involved than the deterministic case in  \cite{tran2019non}, where a new Lyapunov function has been constructed (see \eqref{eq:lyapunov_func}).
}

\vspace{1ex}
\noindent\textbf{Content.}
The rest of this paper is organized as follows.
Section~\ref{sec:math_tools} states our fundamental assumptions and presents some background related to \eqref{eq:primal_form} and \eqref{eq:dual_form}.
Section~\ref{sec:A2_SPD_method} develops our main algorithm, Algorithm~\ref{alg:A1_main}, and establishes its convergence rates in two settings.
This section also investigates two special cases.
Section~\ref{sec_experiment} provides two numerical examples to verify our algorithmic variants and compare them with two other methods.
The proofs of the main results are given in Section~\ref{apdx:sec:appendix1}.

\beforesec
\section{Fundamental Assumptions and Related Background}\label{sec:math_tools}
\aftersec
This section states our fundamental assumption and presents some related background.

\beforesubsec
\subsection{Basic notation and concepts}\label{subsec:basic_notations}
\aftersubsec
We work with finite dimensional spaces $\R^p$, $\R^m$, and $\R^d$, equipped with the standard inner product $\iprods{\cdot,\cdot}$ and Euclidean norm $\norm{\cdot}$.
For any nonempty, closed, and convex set $\Xc$ in $\R^p$, $\relint{\Xc}$ denotes the relative interior of $\Xc$ and $\delta_{\Xc}(\cdot)$ is the indicator of $\Xc$.
For any proper, closed, and convex function $f : \R^p\to\Rext$, $\dom{f}$ denotes its  domain, $f^{\ast}$ is its Fenchel conjugate, $\partial{f}$ denotes its subdifferential \citep{Bauschke2011}. 
We define $\prox_{f}(x) := \textstyle\mathrm{arg}\min_{y}\{ f(y) + (1/2)\norms{y-x}^2 \}$ the proximal operator of $f$.
If $\nabla{f}$ is Lipschitz continuous with a Lipschitz constant $L_f  \geq 0$, i.e., $\norms{\nabla{f}(x) - \nabla{f}(y)}  \leq L_f\norms{x - y}$ for $x, y\in\dom{f}$, then $f$ is called $L_f$-smooth.
If $f(\cdot) - \frac{\mu_f}{2}\norms{\cdot}^2$ is convex for some $\mu_f > 0$, then $f$ is called $\mu_f$-strongly convex with a strong convexity parameter $\mu_f$.
If $\mu_f = 0$, then $f$ is just convex.
We say that $f$ is $M_f$-Lipschitz continuous if $\vert f(x) - f(\hat{x})\vert \leq M_f\norms{x - \hat{x}}$ for all $x, \hat{x}\in\dom{f}$. 
We use $\R_{++}$ to denote the set of positive real numbers, and $[n] := \sets{1, 2, \cdots, n}$ for any positive integer $n$.

Given $K \in \R^{d\times p}$, $K_i$ denotes the $i$-th column block of $K$.
Given $\sigma \in\R^n_{++}$, we define a weighted norm as $\norms{x}_{\sigma} := \big(\sum_{i=1}^n \sigma_i \norms{x_i}^2\big)^{1/2}$ and its dual norm as $\norms{y}_{*} := \big(\sum_{i=1}^n\frac{1}{\sigma_i}\norms{y_i}_{*}^2\big)^{1/2}$.
Let $q \in\R^n_{++}$ be a discrete probability distribution on $[n]$ such that $\sum_{i=1}^nq_i = 1$.
Let $i_k \in [n]$  be a random index such that 
\myeq{eq:ij_distribution}{
\Prob{i_k = i} = q_i.
}
We write $i_k \sim \Up{[n]}$ for sampling a block $i_k$ from $[n]$ based on the distribution $q$.

\beforesubsec
\subsection{Basic assumptions and optimality condition}\label{subsec:assumptions}
\aftersubsec
Our new primal-dual method relies on the following assumptions imposed on  \eqref{eq:primal_form}.

\begin{assumption}\label{as:A0}
The solution set $\Sc^{\star}$ of \eqref{eq:primal_form} is nonempty and the Slater condition $\relint{\dom{f + h}\times\dom{g}\cap\sets{(x, w) \in\R^p\times\R^m : Kx + Bw = b}} \neq \emptyset$ holds.

The functions $f$ and $g$ in \eqref{eq:primal_form} are proper, closed, possibly nonsmooth, and convex on their domain. 
The function $h$ is convex and partially $L_{h,i}$-smooth for all $i\in[n]$, i.e., for any $x\in\R^p$ and $d_i\in\R^{p_i}$ with $i\in [n]$, we have
\myeq{eq:Lh_smooth}{
\norms{\nabla_{x_i}{h}(x + U_id_i) - \nabla_{x_i}{h}(x)}  \leq L_{h,i}\norms{d_i},
}
where $U_i\in\R^{p\times p_i}$ has $p_i$ unit vectors such that $[U_1, U_2, \cdots, U_n]$ forms the identity matrix $\Id$ in $\R^{p\times p}$.
\end{assumption}
Assumption~\ref{as:A0} is often required in primal-dual methods.
Since $\Sc^{\star}$ is nonempty, Assumption~\ref{as:A0} implies strong duality, i.e., $F^{\star} = D^{\star}$, and the solution set $\Yc^{\star}$ of \eqref{eq:dual_form} is also nonempty and bounded.
Let us write $\Sc^{\star} := \Xc^{\star} \times \Wc^{\star}$.

The primal-dual forms \eqref{eq:primal_form} and \eqref{eq:dual_form} can be put into the following minimax form:
\myeq{eq:minimax_form}{
\min_{x\in\R^p, w\in\R^m}\max_{y\in\R^d}\Big\{ \Lc(x, w, y) := F(x, w) + \iprods{Kx + Bw - b, y} \Big\},
}
where $\Lc$ is the Lagrange function associated with \eqref{eq:primal_form} and $y$ is a dual variable or a Lagrange multiplier.
The optimality condition associated with \eqref{eq:primal_form} and its dual form \eqref{eq:dual_form} can be written as follows:
\myeq{eq:opt_cond}{
\hspace{-1ex}
0 \in \partial{f}(x^{\star}) + \nabla{h}(x^{\star}) + K^{\top}y^{\star}, \ \  0 \in \partial{g}(w^{\star}) + B^{\top}y^{\star}, \ \ \text{and} \quad Kx^{\star} + Bw^{\star} -  b = 0.
\hspace{-1ex}
}
Any point $(x^{\star}, w^{\star}, y^{\star})$ satisfying \eqref{eq:opt_cond} is a saddle-point of $\Lc$ in \eqref{eq:minimax_form}, i.e.:
\myeq{eq:saddle_point}{
\hspace{-1ex}
\Lc(x^{\star},  w^{\star}, y) \leq \Lc(x^{\star},  w^{\star}, y^{\star}) \leq \Lc(x, w, y^{\star}), \ \ \forall x \in\dom{f+h},  w\in\dom{g}, y \in \R^d.
\hspace{-1ex}
}
Under Assumption~\ref{as:A0}, $(x^{\star},  w^{\star})$ is  a primal optimal solution of \eqref{eq:primal_form} and $y^{\star}$ is a dual optimal solution of \eqref{eq:dual_form}.

\vspace{1ex}
\noindent\textbf{Primal-dual expected gap.}
To characterize saddle-points of \eqref{eq:minimax_form}, we define
\myeq{eq:gap_func}{
\Gc_{\Zc}(x, w, y) := \sup_{ (\hat{x}, \hat{w}, \hat{y}) \in\Zc}\Exp{ \Lc(x, w, \hat{y}) - \Lc(\hat{x}, \hat{w}, y)},
}
for any nonempty and compact subset $\Zc := \Xc\times \Wc\times \Yc$ in $\R^p\times\R^m\times\R^d$ such that $\Zc \cap \Zc^{\star}\neq\emptyset$, where $\Zc^{\star} = \Xc^{\star} \times \Wc^{\star}\times \Yc^{\star}$.
Here, the expectation is taken overall the randomness generated by the underlying algorithm up to the current iteration for given $(x, w, y)$.
By \eqref{eq:saddle_point}, one can show that $\Gc_{\Zc}(x, w, y) \geq 0$ for any $(x, w, y) \in \R^p \times \R^m\times \R^d$, and if $(x^{\star}, w^{\star}, y^{\star})$ is a saddle-point of \eqref{eq:minimax_form}, then $\Gc_{\Zc} (x^{\star}, w^{\star}, y^{\star}) = 0$.

The function in \eqref{eq:gap_func} has been widely used in convex optimization as well as convex-concave saddle-point problems, see, e.g., \citep{Chambolle2011,tan2018stochastic}.
Note that the expectation in \eqref{eq:gap_func} is inside the supremum instead of outside as the one in  \cite{alacaoglu2019convergence,Chen2013a}.
Hence, we call \eqref{eq:gap_func} a primal-dual expected gap function to distinguish it from \cite{alacaoglu2019convergence,Chen2013a}.
As mentioned in \cite{alacaoglu2019convergence}, there is a technical issue in the proof of \citep{Chambolle2011} and recently in \cite{tan2018stochastic}, leading to an inconsistent conclusion on the gap function guarantee in both papers \cite{Chambolle2011,tan2018stochastic}.

\beforesec
\section{Randomized Block-Coordinate Alternating Primal-Dual Algorithm}\label{sec:A2_SPD_method}
\aftersec
In this section, we develop a unified randomized block-coordinate primal-dual algorithm to solve \eqref{eq:primal_form} and its dual form \eqref{eq:dual_form}.
Then, we investigate its convergence rates.

\beforesubsec
\subsection{The main idea and the full algorithm}\label{subsec:NSPD1}
\aftersubsec
\noindent\mytb{Main idea.} 
Our approach relies on a classical augmented Lagrangian function associated with \eqref{eq:primal_form}, which is defined as follows:
\myeq{eq:aug_Lag}{
\Lc_{\rho}(x, w, y) := f(x) + h(x) + g(w) + \iprods{Kx + Bw -  b, y} + \frac{\rho}{2}\norms{Kx + Bw - b}^2,
}
where $\rho > 0$ is a penalty parameter.
This function will serve as a \textit{\textbf{merit function}} to measure the optimality for both \eqref{eq:primal_form} and its dual form \eqref{eq:dual_form}.

Our central idea can be presented as follows.
\begin{itemize}
\item Firstly, we alternatively minimize $\Lc_{\rho}$ w.r.t. $w$ and $x$.
While the minimization over $w$ is updated in full, the minimization over $x$ is updated by a randomized block-coordinate scheme.
\item  More specifically, the minimization problem over $w$ can be written as
\myeq{eq:w_step1}{
w^{k+1} \in \argmin_{w\in\R^m}\Big\{ g(w) + \iprods{\hat{y}^k, Bw} + \frac{\rho_k}{2}\norms{Bw + K\hat{x}^k - b}^2 \Big\}.
}
\item However, since the minimization problem in $x$ is large-scale, we not only linearize it, but also apply a randomized proximal coordinate gradient method, e.g., in \citep{fercoq2015accelerated} to minimize $\Lc_{\rho_k}(\cdot, w^{k+1}, \hat{y}^k)$.
More concretely, we sample a random block-coordinate $i = i_k$ and update $\tilde{x}_i^k$ by partially linearizing both $\psi_{\rho_k}(x, w^{k+1}, \hat{y}^k) := \iprods{y^k, Kx} + \frac{\rho_k}{2}\norms{Kx + Bw^{k+1} - b}^2$ and $h$ around $\hat{x}^k$ as
\myeqn{
\hspace{-1ex}
\arraycolsep=0.1em
\begin{array}{lcl}
\tilde{x}_{i}^{k+1} & := & {\argmin_{x_i\in\R^{p_i}}}\Big\{f_i(x_i) + \iprods{\nabla_{x_i}{h}(\hat{x}^k) + \nabla_{x_i}{\psi}_{\rho_k}(\hat{x}^k,w^{k+1},\hat{y}^k), x_i - \hat{x}^k_i} \vspace{0.5ex}\\
&& \quad\qquad\qquad + {~}  \frac{\tau_k\sigma_i}{2\tau_0\beta_k}\norms{x_i - \tilde{x}_i^k}^2 \Big\},
\end{array}
\hspace{-3ex}
}
where $\beta_k > 0$ and $\tau_k \in (0, 1)$.
Otherwise, we maintain $\tilde{x}^{k+1}_i := \tilde{x}^k_i$ for $i\neq i_k$.
\rvtext{Here, $\sigma_i$ is a scaling parameter for each block $i$, which will be chosen proportionally to $\norms{K_i}^2$ and $L_{h,i}$ to minimize the overall Lipschitz constants $\bar{L}_{\sigma}$ and $L_{\sigma}^h$ in \eqref{eq:common_quantities}.}
\item  Secondly, we also apply the accelerated steps as in \citep{fercoq2015accelerated} and adaptively update the related parameters $\rho_k$, $\beta_k$, and $\tau_k$ using the ideas in \cite{tran2019non}.
\item Thirdly, we update the dual variable $\hat{y}^k$ instead of fixing it as in \cite{Alacaoglu2017}.
\item Finally, we add an averaging dual step $\bar{y}^k$ to derive dual convergence rates.
\end{itemize}
We now specify each step discussed above to obtain the full algorithm.

\vspace{0.5ex}
\noindent\textbf{The full algorithm.}
Our complete algorithm, called \textit{Randomized block-coordinate alternating primal-dual \textrm{$($PD$)$} algorithm}, is described in detail in Algorithm~\ref{alg:A1_main}.

\begin{algorithm}[ht!]
   \caption{(Randomized Block-Coordinate Alternating Primal-Dual Algorithm)}
   \label{alg:A1_main}
\begin{algorithmic}[1]
   \Statex\label{step:A1-i01}
   \hspace{-4ex}{\bfseries Initialization:} 
   \State\hspace{0ex}Choose $x^0 \in \R^p$ and $\hat{y}^0 \in \R^d$ such that $-B^{\top}\hat{y}^0 \in\dom{g^{*}}$.
   \State\hspace{0ex}Choose $\rho_0 > 0$ or as in Theorem~\ref{th:convergence_rate2}. 
   \State\hspace{0ex}Set $\tilde{x}^0 := x^0$, $\bar{y}^0 := \hat{y}^0$, and $\tau_0 := q_{\min}$, where $q_{\min} := \min_{i\in [n]}q_i > 0$.
   \vspace{0.25ex}   
   \Statex\hspace{-4ex}\textbf{For~}{$k:= 0$ {\bfseries to} $k_{\max}$, \textbf{perform}}
   \vspace{0.5ex}   
   \State\hspace{0ex}\label{step:A1-i1}Update  $\tau_k$, $\rho_k$, $\beta_k$,  and $\eta_k$ as in Theorem~\ref{th:convergence_rate} or Theorem~\ref{th:convergence_rate2}.
   \vspace{0.5ex}   
   \State\hspace{0ex}\label{step:A1-i2a}Update $\hat{x}^{k}  :=  (1-\tau_k)x^k + \tau_k\tilde{x}^k$.
   \vspace{0.5ex}   
   \State\hspace{0ex}\label{step:A1-i2}Update $w^{k+1}$ by solving \eqref{eq:w_step1}.
   \vspace{0.5ex}
   \State\hspace{0ex}\label{step:A1-i7}\rvtext{Update $ \bar{y}^{k+1} := (1-\tau_k)\bar{y}^k + \tau_k\big[\hat{y}^k + \rho_k(K\hat{x}^k + Bw^{k+1} - b)\big]$ (if necessary).}
   \vspace{0.5ex}
   \State\hspace{0ex}\label{step:A1-i3}Sample a block-coordinate $i_k\sim \Up{[n]}$ with the distribution \eqref{eq:ij_distribution}. 
   \vspace{0.5ex}
   \State\hspace{0ex}\label{step:A1-i4}Maintain $\tilde{x}_i^{k+1} := \tilde{x}_i^k$ for all $i\neq i_k$, and for $i = i_k$, update
   \vspace{-1.25ex}
   \myeqn{
   \tilde{x}^{k+1}_i   :=  \prox_{\frac{\tau_0\beta_k}{\tau_k\sigma_i}f_i}\Big(\tilde{x}^k_i - \tfrac{\tau_0\beta_k}{\tau_k\sigma_i}\big(\nabla_{x_i}h(\hat{x}^k) + K^{\top}_i(\hat{y}^k + \rho_k(K\hat{x}^k + Bw^{k+1} -b))\big) \Big).
    \vspace{-1.5ex}
   }
   \State\hspace{0ex}\label{step:A1-i5}Update $x^{k+1}  :=  \hat{x}^k + \frac{\tau_k}{\tau_0}(\tilde{x}^{k+1} - \tilde{x}^k)$.
   \vspace{0.5ex}
   \State\hspace{0ex}\label{step:A1-i6}Update $\hat{y}^{k +1} := \hat{y}^k + \eta_k\big[(Kx^{k +1} + Bw^{k+1} - b) - (1-\tau_k)(Kx^k + Bw^k - b)\big]$.
\Statex\hspace{-4ex}\textbf{EndFor}
\end{algorithmic}
\end{algorithm}

\beforepara
\paragraph*{Per-iteration complexity.}
The main computation of Algorithm~\ref{alg:A1_main} consists of:
\begin{itemize}
\item Step~\ref{step:A1-i2} for updating $w^{k+1}$ requires to solve the subproblem \eqref{eq:w_step1}.
If $B = -\Id$, where $\Id$ is the identity matrix, then $w^{k+1} = \prox_{g/\rho_k}\big(Kx^k - b  + \frac{1}{\rho_k}\hat{y}^k\big)$, which reduces to evaluating one proximal operator of $g$.

\item Step~\ref{step:A1-i7} on $\bar{y}^{k+1}$ is only required if we prove a dual convergence guarantee.
\rvtext{Note that $\bar{y}^{k+1}$ does not depend on $i_k$, and therefore, $\tilde{x}^{k+1}$ and $x^{k+1}$.}

\item Step~\ref{step:A1-i4} only updates one block-coordinate $i_k$ of $\tilde{x}^k$.
This step needs one proximal operation of the component $f_{i_k}$, one partial derivative $\nabla_{x_{i_k}}h$, one $K\hat{x}^k$, one $Bw^{k+1}$, and one $K_{i_k}^{\top}y$.
Note that since ${x}^{k+1}$ is only changed at one block-coordinate $i_{k}$, calculating the product $Kx^{k+1}$ only needs to update $K_{i_{k}}\tilde{x}_{i_{k}}^k$.
\item The dual steps, Steps~\ref{step:A1-i6} and \ref{step:A1-i7}, require to update full vectors in $\R^d$.
\end{itemize}
Currently, we have not specified how to efficiently implement Algorithm~\ref{alg:A1_main}.
We will derive in Subsection~\ref{subsec:impl_variant} an efficient implementation of Algorithm~\ref{alg:A1_main}.

\beforesubsec
\subsection{Convergence guarantees under general convexity}\label{subsec:A1_convergence1}
\aftersubsec
Let us define the following quantities, which will be repeatedly used in the sequel.
\myeq{eq:common_quantities}{
\bar{L}_{\sigma} :=  \max_{i\in [n]}\bigg\{ \frac{\norms{K_i}^2}{\sigma_i} \bigg\}, \quad L_{\sigma}^h := \max_{i \in [n]} \bigg\{ \frac{L_{h,i}}{\sigma_i} \bigg\}, \quad \text{and} \quad \tau_0 := q_{\min} = \min_{i\in [n]}q_i \in (0, 1).
}
We also define $\Fc_k := \sigma(i_0, i_1, \cdots, i_{k-1})$ the $\sigma$-algebra generated by random variables $i_l$ for $l=0,\cdots, k-1$.
We also use a shorthand $\Exps{i_k}{\cdot}$ for the expectation $\Exps{i_k}{\cdot \mid \Fc_{k}}$ conditioned on $\Fc_{k}$, and $\Exp{\cdot}$ for the full expectation on the overall $\sigma$-algebra $\Fc_k$.

We state the first main result for Algorithm~\ref{alg:A1_main} in Theorem~\ref{th:convergence_rate}, whose proof is postponed to Subsection~\ref{apdx:le:convergence_result}.

\begin{theorem}\label{th:convergence_rate}
Suppose that  \eqref{eq:primal_form} satisfies Assumption~\ref{as:A0}, and $\mu_{f_i} = 0$ for all $i \in [n]$.
Let  $\bar{L}_{\sigma}$,  $L_{\sigma}^h$, and $\tau_0$ be given by \eqref{eq:common_quantities} and $\rho_0 > 0$.
Let $\set{(x^k, w^k, \bar{y}^k)}$ be  generated by Algortihm~\ref{alg:A1_main}, where $\tau_k$, $\beta_k$, $\rho_k$, and $\eta_k$ are updated by
\myeq{eq:param_update1}{
\tau_k := \dfrac{\tau_0}{k + 1}, \quad \rho_k := \dfrac{\rho_0\tau_0}{\tau_k}, \quad \beta_k :=  \dfrac{1}{L^h_{\sigma} + 2\bar{L}_{\sigma}\rho_k},  \quad\text{and}\quad \eta_k := \dfrac{\rho_k}{2}.
}
In addition, let $\Gc_{\Zc}$ be defined by \eqref{eq:gap_func}.
Then, the following estimates  hold:
\myeq{eq:residual_bound}{
\hspace{-2ex}
\arraycolsep=0.2em
\left\{\begin{array}{lcll}
\left\vert\Exp{F(x^{k}, w^k) - F^{\star}} \right\vert &\leq& \dfrac{ \bar{\Ec}_0 + \norms{y^{\star}}(2\bar{\Ec}_0/\rho_0)^{1/2}}{\tau_0k + 1 - \tau_0}, &\textrm{$($primal objective residual$)$} \vspace{1ex}\\
\Exp{\norms{Kx^k + Bw^k - b}^2} &\leq& \dfrac{2 \bar{\Ec}_0 }{\rho_0(\tau_0k + 1 - \tau_0)^2}, &\textrm{$($primal feasibility$)$}  \vspace{1ex}\\
\Exp{D^{\star} - D(\bar{y}^k)} &\leq & \dfrac{\bar{\mathcal{F}}_0}{\tau_0k + 1 - \tau_0}, &\textrm{$($dual objective residual$)$}  \vspace{1ex}\\
\Gc_{\Zc}(x^k, w^k, \bar{y}^k) &\leq & \dfrac{F(x^0, w^0) - D(\hat{y}^0) + \bar{R}_{\Zc}^2}{\tau_0k + 1 - \tau_0}, &\textrm{$($primal-dual expected gap$)$}
\end{array}\right.
\hspace{-4ex}
}
where \rvtext{$u^0 := Kx^0 + Bw^0 - b$}, and $\bar{\Ec}_0$, $\bar{\mathcal{F}}_0$, and $\bar{R}_{\Zc}^2$ are respectively defined as
\myeq{eq:Ebar_0}{
\hspace{-3ex}
\arraycolsep=0.1em
\left\{\begin{array}{lcl}
\bar{\Ec}_0 &:=&  F(x^0, w^0) - D(\hat{y}^0)  + \frac{2}{\rho_0}\norms{ y^{\star} - \hat{y}^0}^2 + \rvtext{\frac{1}{\rho_0}\norms{\hat{y}^0}^2} + \frac{\rho_0(2-\tau_0)}{2}\norms{u^0}^2 \vspace{0.5ex}\\
&& + {~} \frac{(L^h_{\sigma} + 2\rho_0\bar{L}_{\sigma})\tau_0}{2}\norms{x^{\star} - x^0}_{\sigma/q}^2, \vspace{1ex}\\
\bar{\mathcal{F}}_0 &:=&  F(x^0, w^0) - D(\hat{y}^0) + \frac{2}{\rho_0}\norms{ y^{\star} - \hat{y}^0}^2 + \rvtext{\frac{1}{\rho_0}\norms{\hat{y}^0}^2} + \frac{\rho_0(2-\tau_0)}{2}\norms{u^0}^2 \vspace{0.5ex}\\
&& + {~} \frac{(L^h_{\sigma} + 2\rho_0\bar{L}_{\sigma})\tau_0}{2}M_0, \vspace{1ex}\\
\bar{R}_{\Zc}^2 &:= &{\!\!\!\!\!\!} {\displaystyle\sup_{(x,y)\in\Xc\times\Yc}} {\!} \Big\{ \tfrac{(L^h_{\sigma} + 2\rho_0\bar{L}_{\sigma})\tau_0}{2}\norms{x - x^0}_{\sigma/q}^2 +  \tfrac{2}{\rho_0}\norms{y - \hat{y}^0}^2 \Big\} + \rvtext{\frac{1}{\rho_0}\norms{\hat{y}^0}^2} + \frac{\rho_0(2-\tau_0)}{2}\norms{u^0}^2.
\end{array}\right.
\hspace{-8ex}
}
Here, we assume that the Fenchel conjugate $\phi^{*}$ of $\phi = f + h$ is $M_{\phi^{*}}$-Lipschitz continuous and $M_0 := \sup\sets{\norms{x - x^0}^2 : \norms{x} \leq M_{\phi^{*}}}$ to obtain the dual objective residual bound.
\end{theorem}

\rvtext{
\begin{remark}[\textit{\textbf{The fininiteness of $D(\bar{y}^k)$}}]\label{re:finite_Dy}
Note that, by \eqref{eq:dual_form}, we have
\begin{equation*}
D(y) = \min_{x,w}\left\{\phi(x) + g(w) + \iprod{Kx + Bw -b, y}\right\} = -\phi^{\ast}(-K^{\top}y) - g^{\ast}(-B^{\top}y) - \iprod{b, y}.
\end{equation*}
Hence, we have $\dom{D} = \sets{y \in \R^d : -K^{\top}y \in \dom{\phi^{\ast}}, \ -B^{\top}y \in \dom{g^{\ast}}}$.
We show in the proof of Theorem \ref{th:convergence_rate} that $\{\bar{y}^k\}$ in Algorithm \ref{alg:A1_main} always belongs to $\dom{D}$ as long as $-B^{\top}\bar{y}^0 \in \dom{g^{\ast}}$ and $\phi^{\ast}$ is $M_{\phi^{\ast}}$-Lipschitz continuous.
Consequently, since $D$ is concave and proper, we conclude that $D(\bar{y}^k)$ is finite.
\end{remark}
}

\begin{remark}\label{re:optimal_param}
If we choose $q_i := \frac{1}{n}$ for all $i \in [n]$, then $\tau_0 = \frac{1}{n}$, and the convergence rate in Theorem~\ref{th:convergence_rate} is $\BigO{\frac{n}{k}}$, which matches the rate in \cite{Alacaoglu2017}. If $n=1$, then this rate is optimal up to a constant factor as discussed in \cite{tran2019non}.

If we choose
\myeqn{
\sigma_i  = q_i = \frac{L_{h,i}  + \rho_0\norms{K_i}^2}{\sum_{i=1}^{n}(L_{h,i} + \rho_0\norms{K_i}^2)}, \quad \forall i\in [n],
}
then, Algorithm~\ref{alg:A1_main} takes into account the Lipschitz constant $L_{h,i}$ of $h$ and $\norms{K_i}^2$ of each block $i$.
This is expected to improve the performance of Algorithm~\ref{alg:A1_main} when the input data represented in $K$ is not normalized, and the partial Lipschitz constants $L_{h,i}$ is really different between block-coordinates $x_i$ for some $i\in [n$].
\end{remark}

\beforesubsec
\subsection{Convergence guarantees under strong convexity}\label{subsec:A1_convergence2}
\aftersubsec
If either $f$ or $h$ in \eqref{eq:primal_form} is strongly convex, then we can boost Algorithm~\ref{alg:A1_main} up to  $\BigOc{1/k^2}$ convergence rate.
The following theorem states this acceleration when $f$ is strongly convex, whose proof is deferred to Subsection~\ref{apdx:th:convergence_rate2}.

\begin{theorem}\label{th:convergence_rate2}
Suppose that \eqref{eq:primal_form} satisfies Assumption~\ref{as:A0} and $f$ is strongly convex, i.e., $\mu_{f_i} > 0$ for all $i \in [n]$, but $h$ and $g$ are not necessarily strongly convex.
Let $\tau_0$, $\bar{L}_{\sigma}$, and $L^h_{\sigma}$ be given by \eqref{eq:common_quantities} and $\set{(x^k, w^k, \bar{y}^k)}$ be generated by Algorithm~\ref{alg:A1_main}, and $\tau_k$,  $\beta_k$, $\rho_k$, and $\eta_k$ are updated by 
\myeqn{ 
\tau_k := \frac{\tau_{k-1}\big[(\tau_{k-1}^2+4)^{1/2} - \tau_{k-1}\big]}{2}, \ \ \rho_k := \dfrac{\rho_{k-1}}{1-\tau_k}, \ \  \beta_k := \dfrac{1}{L^h_{\sigma} + 2\bar{L}_{\sigma}\rho_k}, \ \ \text{and} \ \ \eta_k := \dfrac{\rho_k}{2},
}
where $\rho_0$ is chosen such that $0 < \rho_0 \leq \frac{1}{4\bar{L}_{\sigma}}\min\big\{ \frac{\mu_{f_i}}{\sigma_i} : i \in [n] \big\}$.
Then, we have
\myeq{eq:key_estimate3_scvx}{
\arraycolsep=0.2em
\hspace{-2ex}
\left\{\begin{array}{lcll}
\left\vert\Exp{F(x^{k}, w^k) - F^{\star}} \right\vert &\leq& \dfrac{ 4\big[\tilde{\Ec}_0 + \norms{y^{\star}}(2\tilde{\Ec}_0/\rho_0)^{1/2}\big]}{(\tau_0k + 2)^2}, &\textrm{$($primal objective residual$)$} \vspace{1ex}\\
\Exp{\norms{Kx^k + Bw^k - b}^2} &\leq& \dfrac{8 \tilde{\Ec}_0 }{\rho_0(\tau_0k + 2)^4}, &\textrm{$($primal feasibility$)$}  \vspace{1ex}\\
\Exp{D^{\star} - D(\bar{y}^k)} &\leq & \dfrac{4\tilde{\Fc}_0}{(\tau_0k + 2)^2}, &\textrm{$($dual objective residual$)$}  \vspace{1ex}\\
\Gc_{\Zc}(x^k, w^k, \bar{y}^k) &\leq & \dfrac{F(x^0, w^0) - D(\hat{y}^0) + \tilde{R}_{\Zc}^2}{(\tau_0k + 2)^2}, &\textrm{$($primal-dual expected gap$)$},
\end{array}\right.
\hspace{-4ex}
}
where \rvtext{$u^0 := Kx^0 + Bw^0 - b$}, and $\tilde{\Ec}_0$, $\tilde{\Fc}_0$, and $\tilde{R}_{\Zc}^2$ are respectively defined as
\myeqn{
\arraycolsep=0.1em
\hspace{-0.5ex}
\left\{\begin{array}{lcl}
\tilde{\Ec}_0 &:=&  F(x^0, w^0) - D(\hat{y}^0)   
+ \sum_{i=1}^n\tfrac{\tau_0}{2q_i}\big[(L^h_{\sigma} + 2\rho_0\bar{L}_{\sigma})\sigma_i  + \mu_{f_i} \big] \norms{x_i^{\star} - x^0_i}^2 \vspace{1ex}\\
&&  + {~} \frac{2}{\rho_0}\norms{ y^{\star} - \hat{y}^0}^2 + \rvtext{\frac{1}{\rho_0}\norms{\hat{y}^0}^2} + \frac{\rho_0(2-\tau_0)}{2}\norms{u^0}^2, \vspace{1ex}\\
\tilde{\Fc}_0 &:=&  F(x^0, w^0) - D(\hat{y}^0)   
+ \sum_{i=1}^n\tfrac{\tau_0}{2q_i}\big[(L^h_{\sigma} + 2\rho_0\bar{L}_{\sigma})\sigma_i  + \mu_{f_i} \big]M_0 \vspace{1ex}\\
&& + {~} \frac{2}{\rho_0}\norms{ y^{\star} - \hat{y}^0}^2 + \rvtext{\frac{1}{\rho_0}\norms{\hat{y}^0}^2 + \frac{\rho_0(2-\tau_0)}{2}\norms{u^0}^2}, \vspace{1ex}\\
\tilde{R}_{\Zc}^2 &:= & {\!\!\!} {\displaystyle\sup_{(x,y)\in\Xc\times\Yc}}  \left\{ \sum_{i=1}^n\tfrac{\tau_0}{2q_i}\big[(L^h_{\sigma} + 2\rho_0\bar{L}_{\sigma})\sigma_i  + \mu_{f_i} \big] \norms{x_i - x^0_i}^2 + \tfrac{2}{\rho_0}\norms{y - \hat{y}^0}^2 \right\} \vspace{1ex}\\
&& \rvtext{+ {~} \frac{1}{\rho_0}\norms{\hat{y}^0}^2 + \frac{\rho_0(2-\tau_0)}{2}\norms{u^0}^2.}
\end{array}\right.
\hspace{-3ex}
}
Here, we assume that the Fenchel conjugate $\phi^{*}$ of $\phi = f + h$ is $M_{\phi^{*}}$-Lipschitz continuous and $M_0 := \sup\sets{\norms{x - x^0}^2 : \norms{x} \leq M_{\phi^{*}}}$ to obtain the dual objective residual bound.
\end{theorem}

Note that we can choose $\hat{y}^0 := 0$ in Theorem~\ref{th:convergence_rate} and Theorem~\ref{th:convergence_rate2} to simplify our corresponding convergence bounds as long as $0 \in \dom{g^{*}}$.
Similar to Theorem~\ref{th:convergence_rate}, if we choose $q_i = \frac{1}{n}$ for all $i\in [n]$, then the convergence rates of Theorem~\ref{th:convergence_rate2} is $\BigO{n^2/k^2}$, which is optimal (up to a constant factor) by assuming that $n=1$ as shown in \cite{tran2019non}.
We can also choose another $q_i$ as in Remark~\ref{re:optimal_param}.

\vspace{1ex}
\noindent\textbf{Handling strong convexity of $h$.}
In Theorem~\ref{th:convergence_rate2}, we assume that $f$ is strongly convex, but $h$ is not necessarily strongly convex.
However, if $f$ is just convex, but $h$ is $\mu_h$-strongly convex with $\mu_h > 0$, then we can process as follows:
\begin{itemize}
\item Replace $h$ by $\tilde{h}(x) := h(x) - \frac{\mu_h}{2}\norms{x}^2$, which is only convex.
Moreover, $\tilde{h}$ is also $L_{\tilde{h},i}$-smooth w.r.t. $x_i$ with $L_{\tilde{h},i} = L_{h,i} - \mu_h$, and $\nabla_{x_i}\tilde{h}(x) = \nabla_{x_i}h(x) - \mu_hx_i$.
\item Replace $f$ by $\tilde{f}(x) := f(x) + \frac{\mu_h}{2}\norms{x}^2$, which is $\mu_h$-strongly convex.
In this case, $\tilde{f}_i(x_i) = f_i(x_i) + \frac{\mu_h}{2}\norms{x_i}^2$, which is also $\mu_h$-strongly convex for all $i\in [n]$.
Moreover, we have $\prox_{\gamma \tilde{f}_i}(x_i) = \prox_{\gamma f_i/(1 + \gamma\mu_h)}(x_i/(1+\gamma \mu_h))$.
\end{itemize}
With these modifications, we apply Algorithm~\ref{alg:A1_main} to $F(x,w) = \tilde{f}(x) + \tilde{h}(x) + g(w)$, and the convergence guarantees in Theorem~\ref{th:convergence_rate2} still hold for this case.

\beforesubsec
\subsection{Two special variants of Algorithm~\ref{alg:A1_main}}\label{subsec:special_variants}
\aftersubsec
Let us consider two special variants of Algorithm~\ref{alg:A1_main} which cover several existing works. 

\beforesubsubsec
\subsubsection{Nonsmooth constrained convex problem}\label{subsubsec:special_case1}
\aftersubsubsec
If $g = 0$ and $B = 0$, then problem \eqref{eq:primal_form} reduces to 
\myeq{eq:primal_form1}{
\begin{array}{l}
{\displaystyle\min_{x\in\R^p}} \Big\{ F(x) := \sum_{i=1}^nf_i(x_i) + h(x) \quad \textrm{s.t.} \quad Kx = b \Big\}.
\end{array}
}
In this case, the main steps of Algorithm~\ref{alg:A1_main} for solving \eqref{eq:primal_form1} can be written as follows:
\myeq{eq:variant1}{
\hspace{-2ex}
\arraycolsep=0.1em
\left\{\begin{array}{lcl}
\tilde{x}^{k+1}_i   & := &  \left\{\begin{array}{ll} 
&\prox_{\frac{\tau_0\beta_k}{\tau_k\sigma_i}f_i}\Big(\tilde{x}^k_i - \tfrac{\tau_0\beta_k}{\tau_k\sigma_i}\big(\nabla_{x_i}h(\hat{x}^k) + K^{\top}_i(\hat{y}^k + \rho_k(K\hat{x}^k -b))\big) \Big), ~ \text{if $i=i_k$},\vspace{1ex}\\
& \tilde{x}^k_i, ~~~\text{otherwise},
\end{array}\right.\vspace{1ex}\\
\hat{y}^{k +1} & := & \hat{y}^k + \eta_k\big[(Kx^{k +1} - b) - (1-\tau_k)(Kx^k  - b)\big]. \vspace{1ex}\\
\end{array}\right.
\hspace{-4ex}
}
Other steps remain the same as in Algorithm~\ref{alg:A1_main}, except that we remove all $w^k$.

Note that \eqref{eq:primal_form1} covers the model in \citep{tan2018stochastic} as a special case by appropriately reformulating it into \eqref{eq:primal_form1}.
More specifically, to process $f(u) + g(Mu)$, we write $f(u) + g(v)$ subject to $Mu - v = 0$.
If we define $x := (u, v)$, then we can transform the model in \citep{tan2018stochastic} into \eqref{eq:primal_form1}.
In this case, the variant \eqref{eq:variant1} of Algorithm~\ref{alg:A1_main} has fully randomized block-coordinate updates overall the primal variable $x$.
To align with \citep[Algorithm 1]{tan2018stochastic}, we can choose a pair of random block-coordinates $(i_k, j_k)$ for $i_k \in [n_1]$ and $j_k \in [n_2]$, where $n_1$ is the number of blocks in $u$ and $n_2$ is the number of blocks in $v$. Here, we update blocks $i_k$ and $j_k$ simultaneously.
In terms of convergence rates, \citep[Algorithm 1]{tan2018stochastic} only achieves $\BigO{1/\sqrt{k}}$ rate, while Algorithm~\ref{alg:A1_main} has a much better rate, which is $\BigO{1/k}$ as shown in Theorem~\ref{th:convergence_rate}. 
Moreover, our convergence rate is non-ergodic (i.e., on the last primal iterate $x^k$) as opposed to averaging sequences in \cite{tan2018stochastic}.

\beforesubsubsec
\subsubsection{Nonsmooth composite convex minimization}\label{subsubsec:special_case2}
\aftersubsubsec
If $b = 0$ and $B = -\Id$, then \eqref{eq:primal_form} reduces to the following setting:
\myeq{eq:primal_form2}{
\begin{array}{l}
{\displaystyle\min_{x\in\R^p}}  \Big\{ F(x) := \sum_{i=1}^nf_i(x_i) + h(x) + g(Kx) \Big\}.
\end{array}
}
In this case, from \eqref{eq:w_step1}, we have $w^{k+1} := \prox_{g/\rho_k}\big(\frac{1}{\rho_k}\hat{y}^k + K\hat{x}^k \big)$.
Using Moreau's identity, we have $w^{k+1} = \frac{1}{\rho_k}(\hat{y}^k + \rho_kK\hat{x}^k - y^{k+1})$ with $y^{k+1} := \prox_{\rho_kg^{*}}\big(\hat{y}^k + \rho_kK\hat{x}^k\big)$.
Utilizing these relations, the main steps of Algorithm~\ref{alg:A1_main} for solving \eqref{eq:primal_form2} can be written as
\myeq{eq:variant2}{
\arraycolsep=0.1em
\left\{\begin{array}{lcl}
y^{k+1} &:= & \prox_{\rho_kg^{*}}\big(\hat{y}^k + \rho_kK\hat{x}^k \big), \vspace{1ex}\\
\tilde{x}^{k+1}_i   & := &  \left\{\begin{array}{lll} 
&\prox_{\frac{\tau_0\beta_k}{\tau_k\sigma_i}f_i}\Big(\tilde{x}^k_i - \tfrac{\tau_0\beta_k}{\tau_k\sigma_i}\big(\nabla_{x_i}h(\hat{x}^k) + K^{\top}_iy^{k+1}\big) \Big), &~~ \text{if $i=i_k$},\vspace{1ex}\\
& \tilde{x}^k_i &~~\text{otherwise},
\end{array}\right. 
\end{array}\right.
}
This step is similar to the main step of SPDHG and other existing primal-dual methods, see, e.g., \cite{Alacaoglu2017,chambolle2017stochastic}.
The convergence of this variant can be derived from Theorems~\ref{th:convergence_rate} and \ref{th:convergence_rate2} combining with the results in \cite{tran2019non}. 
But we omit the details here to avoid overloading this paper.
More specifically,  the primal convergence is given as follows:
\myeqn{
0 \leq \mathbb{E}\big[ F(x^k) - F^{\star}\big]  \leq \left\{\begin{array}{ll}
\frac{ \bar{\Ec}_0 + (\norms{y^{\star}} + M_g)\sqrt{2\bar{\Ec}_0/\rho_0}}{\tau_0k + 1 - \tau_0} &\text{if $\min_i\mu_{f_i} = 0$}, \vspace{1ex}\\
\frac{4\big[ \tilde{\Ec}_0 + (\norms{y^{\star}} + M_g)\sqrt{2\tilde{\Ec}_0/\rho_0}\big]}{(\tau_0k + 2)^2} &\text{if $\min_i\mu_{f_i} > 0$},
\end{array}\right.
}
where $g$ is $M_g$-Lipschitz continuous, and $\bar{\Ec}_0$ and $\tilde{\Ec}_0$ are respectively defined in Theorems~\ref{th:convergence_rate} and \ref{th:convergence_rate2}, but using $F(x^0)$ instead of $F(x^0, w^0)$.
We can also derive convergence rates on the dual objective residual and primal-dual expected gap for \eqref{eq:primal_form2} (see \cite{tran2019non}).

The setting \eqref{eq:primal_form2} covers the models studied, e.g., in \cite{Alacaoglu2017,chambolle2017stochastic,tran2019non,zhang2017stochastic}.
However, the variant \eqref{eq:variant2} is still different from \cite[Algorithm~1]{Alacaoglu2017}, where it has dual updates $\hat{y}^k$ and $\bar{y}^k$.
In \cite[Algorithm 1]{Alacaoglu2017}, $\hat{y}^k$ is fixed at $\dot{y}$ without any dual update, making it less flexible to monitor the dual progress.
Moreover, \cite{Alacaoglu2017} only considers the general convex case and primal convergences, while we consider both the general case and the strongly convex case, and also prove three convergence criteria.
Compared to \cite{alacaoglu2019convergence,chambolle2017stochastic,zhang2017stochastic}, the variant \eqref{eq:variant2} of Algorithm~\ref{alg:A1_main} relies on a Nesterov's accelerated scheme, which has optimal convergence rates on the last iterate $x^k$ as opposed to ergodic sequences as in these works.
Note also that, by eliminating $w^k$ in the update of $\hat{y}^k$, then $\hat{y}^k$ will depend on three consecutive iterates at the iterations $k+1$, $k$, and $k-1$ as discussed in \cite{tran2019non}.

\beforesubsec
\subsection{Efficient implementation of Algorithm~\ref{alg:A1_main}}\label{subsec:impl_variant}
\aftersubsec
In order to efficiently implement Algorithm~\ref{alg:A1_main}, we introduce three intermediate vectors $\tilde{u}^k := K\tilde{x}^k$, $u^k := Kx^k$, and $v^k := Bw^k - b$ in $\R^d$ to store matrix-vector products.
Then, we have $K\hat{x}^k = (1-\tau_k)u^k + \tau_k\tilde{u}^k$.
Hence, for $i=i_k$, we obtain 
\myeqn{
\tilde{x}_{i}^{k+1} := \prox_{\frac{\tau_0\beta_k}{\tau_k\sigma_i}f_i}\big(\tilde{x}^k_i - \tfrac{\tau_0\beta_k}{\tau_k\sigma_i}\big(\nabla_{x_i}h(\hat{x}^k) + K^{\top}_iy^{k+1}\big) \big),
}
where  $y^{k+1} :=  (1-\tau_k)u^k + \tau_k\tilde{u}^k + v^{k+1}$.

Next, using $\tilde{u}^k$, $u^k$, $v^k$, and $y^{k+1}$ above, we can update other steps as follows:
\begin{itemize}
\item Update $u^{k+1} := (1-\tau_k)u^k + \tau_k\tilde{u}^k + \frac{\tau_k}{\tau_0}(\tilde{u}^{k+1} - \tilde{u}^k)$.
\item Update $\hat{y}^{k +1} := \hat{y}^k + \eta_k\big[u^{k+1} + v^{k+1} - (1-\tau_k)(u^k + v^k)\big]$.
\item Update the averaging dual vector: $\bar{y}^{k+1} := (1-\tau_k)\bar{y}^k + \tau_ky^{k+1}$.
\end{itemize}
By exploiting the tricks in \cite{Alacaoglu2017,fercoq2015accelerated}, we can remove the full vector operations in $\R^p$, and only perform the updates on each block $i_k \in [n]$ at each iteration $k$.
We omit this derivation and refer to  \cite{Alacaoglu2017,fercoq2015accelerated} for more details.

\beforesec
\section{Numerical Experiments}\label{sec_experiment}
\aftersec
In this section, we provide two numerical examples to verify our theoretical results and compare Algorithm~\ref{alg:A1_main} with some recent existing methods.
Our first aim is to verify the theoretical convergence rates of Algorithm~\ref{alg:A1_main} stated in Theorems~\ref{th:convergence_rate} and \ref{th:convergence_rate2}.
Then, we compare Algorithm~\ref{alg:A1_main} with two other candidates: SPDHG \citep{chambolle2017stochastic} and PDHG \citep{Chambolle2011} on two common examples. 
Note that both Algorithm \ref{alg:A1_main} and SPDHG essentially have the same per-iteration complexity and convergence rates, but the last-iterate vs. ergodic.

We implement our methods in Python and adapt the code of SPDHG and PDHG from \url{https://github.com/mehrhardt/spdhg}. 
Our experiments were run on a Linux desktop with 3.6GHz Intel Core i7-7700 and 16Gb memory.

\vspace{0.75ex}
\noindent\textbf{Parameter selection strategies.}
For Algorithm~\ref{alg:A1_main}, we search the initial value $\rho_0$ in the range $[1/\norms{K}, 0.1]$ for each dataset, and other parameters are updated by exactly following the update rules in Theorem~\ref{th:convergence_rate} and  Theorem~\ref{th:convergence_rate2}, respectively without any tuning.
For PDHG and SPDHG, we finely tune their step-sizes $\tau$ and  $\sigma$  in the range $[1/\norms{K}, 0.1]$. 
We also tune the extrapolation parameter $\theta$ in the range $[1, d]$ for each dataset, where $d$ is the number of rows of matrix $K$.
We pick the best values of the parameters after tuning for each algorithm to perform our experiments.
By default, the number of block-coordinates is chosen as $n= 32$ in all algorithms.
However, we also use other choices to examine the performance of Algorithm~\ref{alg:A1_main}.

\beforesubsec
\subsection{Support vector machine}\label{subsec_svm}
\aftersubsec
Given a training set of $m$ examples $\{(a_i, b_i)\}_{i=1}^m$, $a_i \in \R^p$ and class labels $b_i \in \{-1, +1\}$, the soft margin SVM problem (without bias) is defined as
\myeq{eq:svm_problem}{
\min_{x \in\R^p}\Big\{F(x) :=  \tfrac{1}{m}\sum_{i = 1}^{m} \max\set{0, 1 - b_i\iprod{a_i, x}} +  \tfrac{\lambda}{2}\norms{x}^2\Big\}.
}
Let us define $g(w) := \frac{1}{m}\sum_{i=1}^m\max\set{0, 1 - w_i}$, $f(x) := \frac{\lambda}{2}\norms{x}^2$, $h(x) := 0$, and using a linear constraint $Ax - w = 0$, where $b_ia_i$ is the $i$-th row of $A$.
Then, \eqref{eq:svm_problem} can be cast into \eqref{eq:primal_form}.
To perform our tests, we use several real datasets from LIBSVM \citep{CC01a}.

\begin{figure*}[ht!]
\begin{center}
\includegraphics[width=0.48\textwidth]{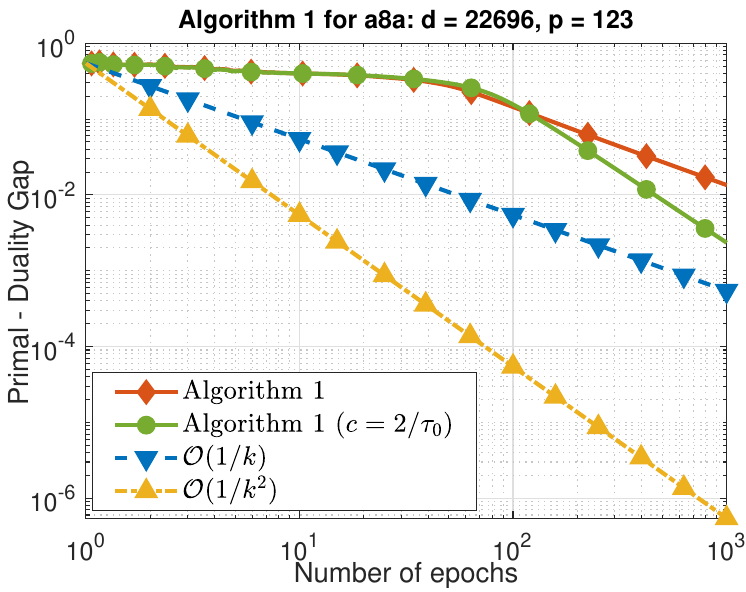}
\includegraphics[width=0.48\textwidth]{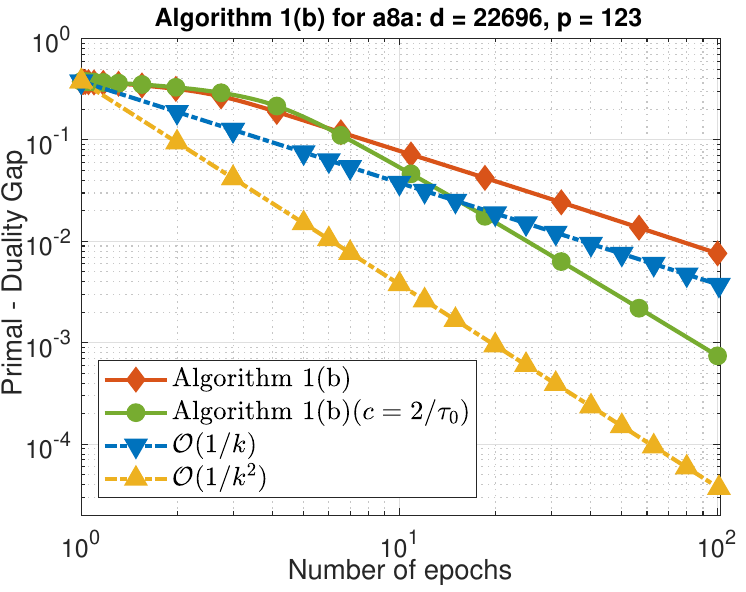}
\includegraphics[width=0.48\textwidth]{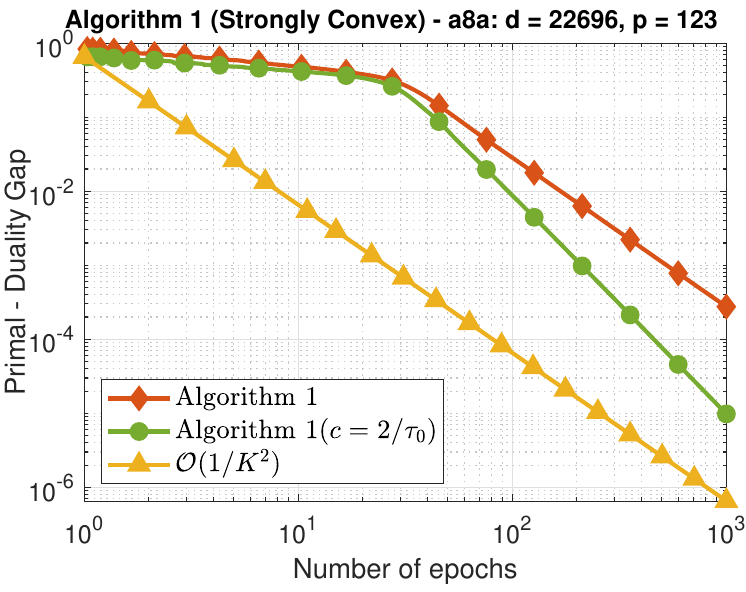}
\includegraphics[width=0.48\textwidth]{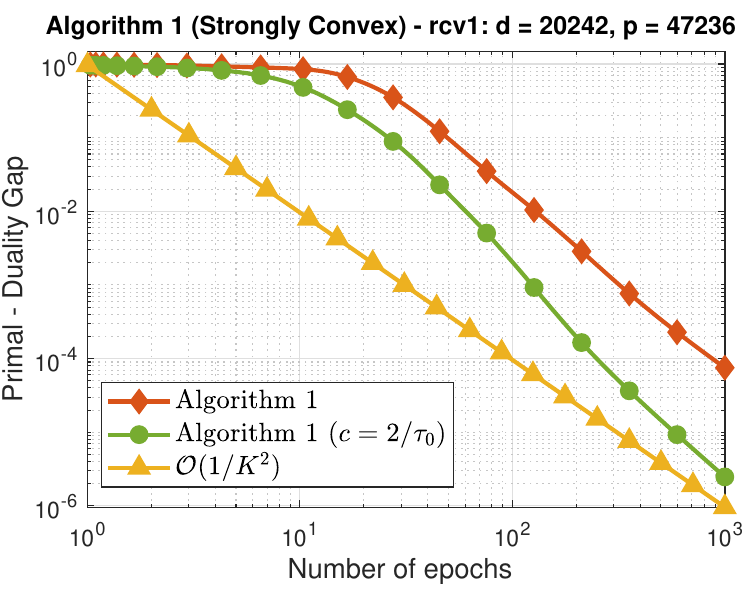}
\hspace{-2ex}
\caption{Convergence rates of Algorithm \ref{alg:A1_main}(Theorems~\ref{th:convergence_rate} and \ref{th:convergence_rate2}) and its variant, Algorithm~\ref{alg:A1_main}(b), and using a  modified rule of $\tau_k$ for solving \eqref{eq:svm_problem} on the \mytb{a8a} and \mytb{rcv1} datasets.
}
\label{fig:compare_c}
\end{center}
\vspace{-2ex}
\end{figure*}

\vspace{1ex}
\noindent\textbf{Experiment 1: Theoretical rate illustration.}
We first illustrate the $\BigOc{1/k}$ convergence rate of Algorithms \ref{alg:A1_main} in Theorem \ref{th:convergence_rate} to solve \eqref{eq:svm_problem} using the \mytb{a8a} dataset in LIBSVM \citep{CC01a}.
Figure~\ref{fig:compare_c} (the top-left plot) shows the convergence behavior of Algorithm~\ref{alg:A1_main} when $\mu_f = 0$ (corresponding to Theorem~\ref{th:convergence_rate}) on the duality gap $F(x^k) - D(\bar{y}^k)$ (which has the same rate as $F(x^k) - F^{\star}$ and $D^{\star} - D(\bar{y}^k)$).
Here, we also modify the update of $\tau_k$ in Theorem~\ref{th:convergence_rate} by $\tau_k := \frac{c\tau_0}{k+c}$ for $c :=2/\tau_0$ to observe faster convergence rates as in deterministic algorithms, see, e.g., \cite{tran2019non}. The plot of this variant is in green.

It is interesting to see that without tuning $\rho_0$, Algorithm \ref{alg:A1_main} converges with $\BigOs{1/k}$-rate if $\mu_f = 0$ as stated by Theorem~\ref{th:convergence_rate}.
If we modify $\tau_k$ with $c := 2/\tau_0$ as mentioned, a faster rate is observed in the green curve.

Now, we reformulate \eqref{eq:svm_problem} into \eqref{eq:primal_form1} in order to test the variant \eqref{eq:variant1} of Algorithm~\ref{alg:A1_main}, denoted by Algorithm~\ref{alg:A1_main}(b).
The top-right plot of Figure~\ref{fig:compare_c} shows the performance of this variant. We still obtain similar convergence rates as shown in the first test.

Finally, we test the $\BigOc{1/k^2}$ rates for Algorithm~\ref{alg:A1_main} when $\mu_f := \lambda > 0$ in \eqref{eq:svm_problem} on the \texttt{a1a} and \texttt{rcv1} datasets.
The result is shown in the bottom-left and bottom-right plots of Figure~\ref{fig:compare_c}.
With the update rules of parameters as in Theorem~\ref{th:convergence_rate2}, we obtain $\BigOc{1/k^2}$ rate as theoretically stated.
This actual rate can be boosted faster than $\BigO{1/k^2}$ if we choose $\tau_k := \frac{c\tau_0}{k+c}$ for $c :=2/\tau_0$ (green curves).

\begin{figure*}[ht!]
\begin{center}
\hspace{-2ex}
\includegraphics[width=0.33\textwidth]{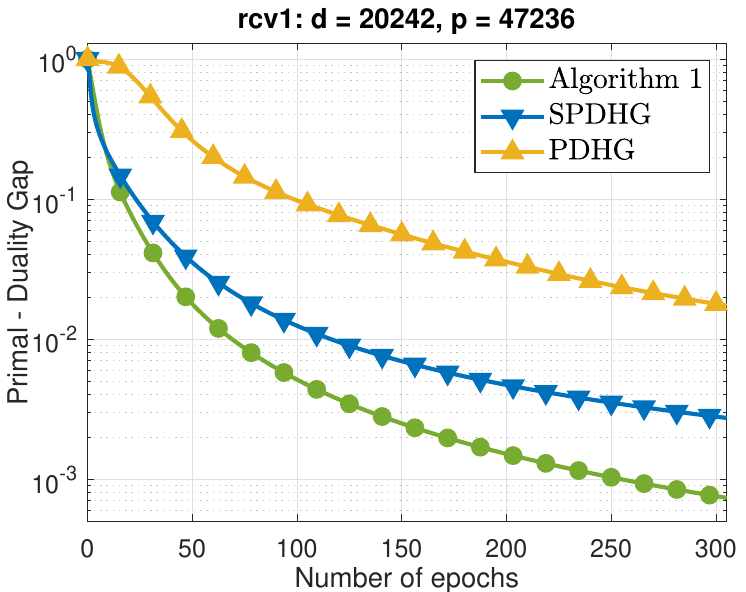}
\includegraphics[width=0.33\textwidth]{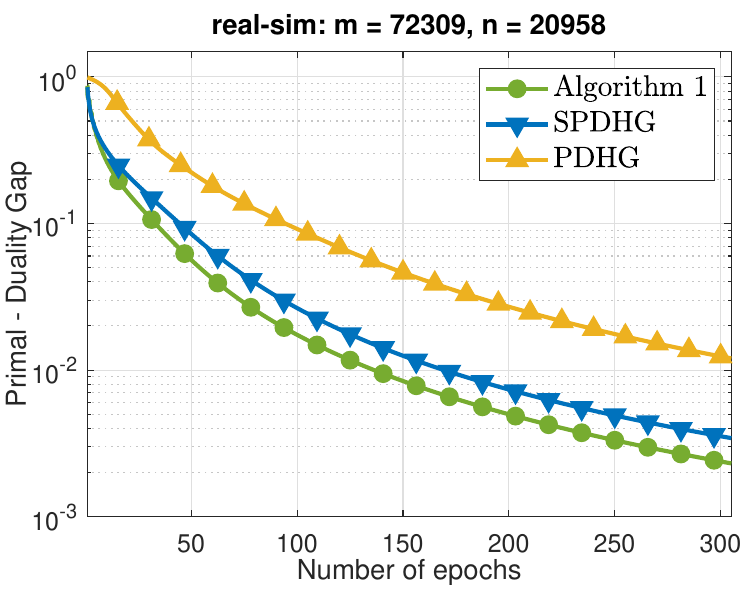}
\includegraphics[width=0.33\textwidth]{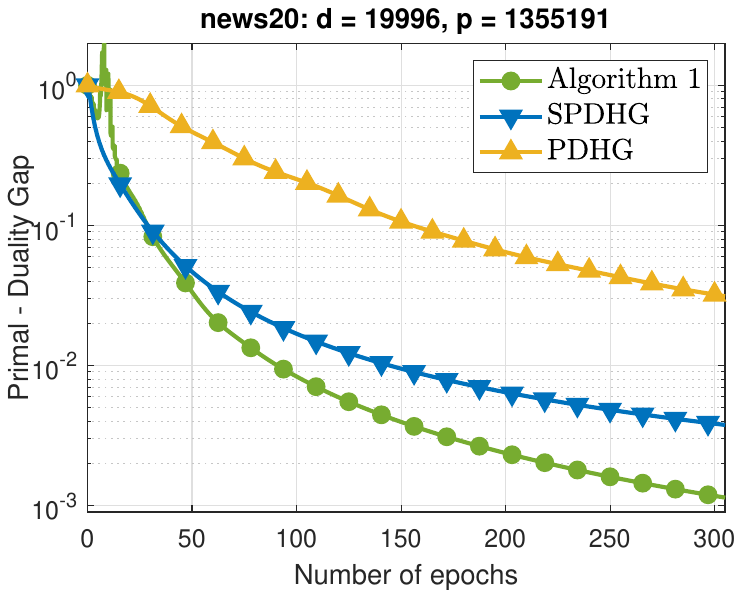}
\hspace{-2ex}
\caption{Comparison of three algorithms, including Algorithm~\ref{alg:A1_main}, for solving \eqref{eq:svm_problem} on $3$ different datasets.}
\label{fig:svm}
\end{center}
\vspace{-3ex}
\end{figure*}

\vspace{1ex}
\noindent\textbf{Experiment 2: Comparison with PDHG and SPDHG.}
We apply Algorithm \ref{alg:A1_main} to solve \eqref{eq:svm_problem} and compare it with SPDHG \citep{chambolle2017stochastic} and PDHG \citep{Chambolle2011,Esser2010}.
We observe that SPDHG is almost identical to SPDC in \citep{zhang2017stochastic} except for assumptions. 
We only choose the variant \eqref{eq:variant2} of Algorithm \ref{alg:A1_main} since it has almost the same per-iteration complexity as SPDHG.
However, we do not take into account the strong convexity of $f$ in this test.
We have tuned these algorithms to obtain the best parameter setting for each dataset.
We test all these algorithms on three different datasets in LIBSVM: \mytb{rcv1}, \mytb{real-sim}, and \mytb{news20} and set $\lambda$ to $10^{-4}$.
The performance is shown in Figure \ref{fig:svm}, where the duality gap $F(x^k) - D(\bar{y}^k)$ is used to measure the performance.  

From Figure \ref{fig:svm}, we can see that Algorithm~\ref{alg:A1_main} gives better convergence behavior than SPDHG in all the datasets.
As usual, stochastic variants such as Algorithm~\ref{alg:A1_main} and SPDHG outperform the deterministic variant, PDHG.
In Figure \ref{fig:svm}, the stochastic algorithms are implemented by separating the whole dimensions into $n=32$ blocks and updating one block at each iteration. 

\begin{figure*}[ht!]
\hspace{-2ex}
\begin{center}
\includegraphics[width=0.328\textwidth]{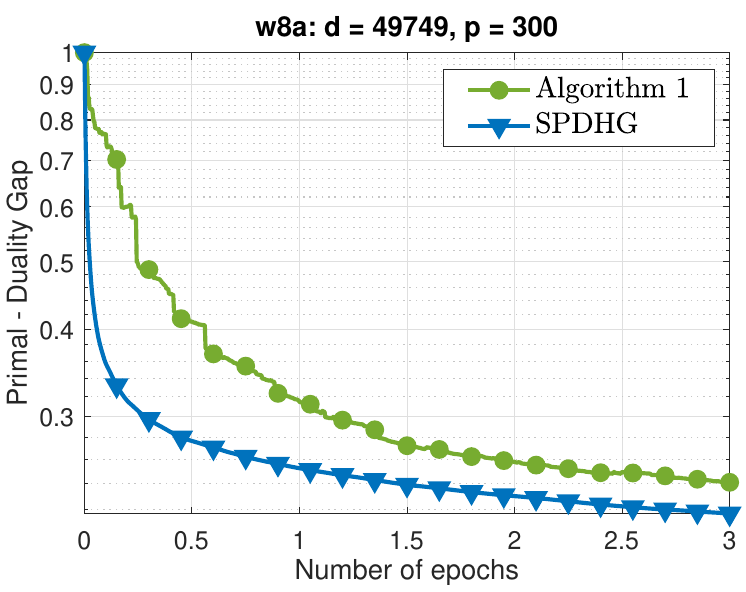}
\includegraphics[width=0.328\textwidth]{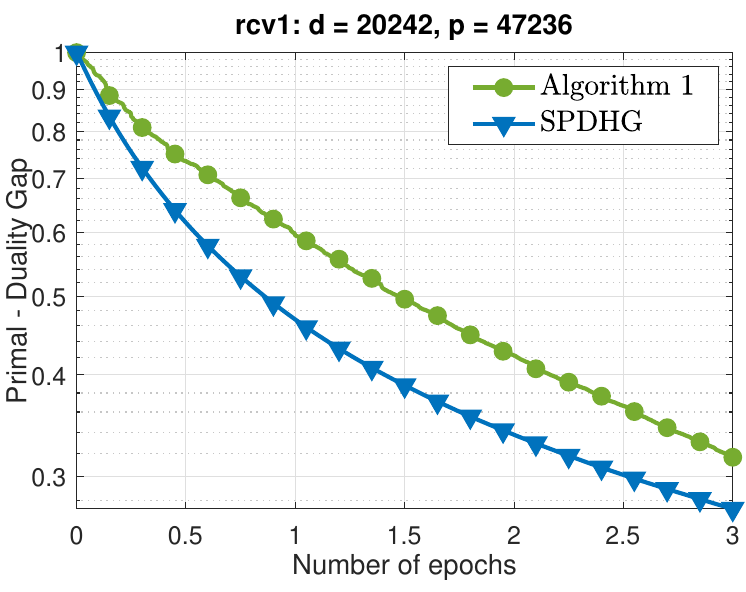}
\includegraphics[width=0.328\textwidth]{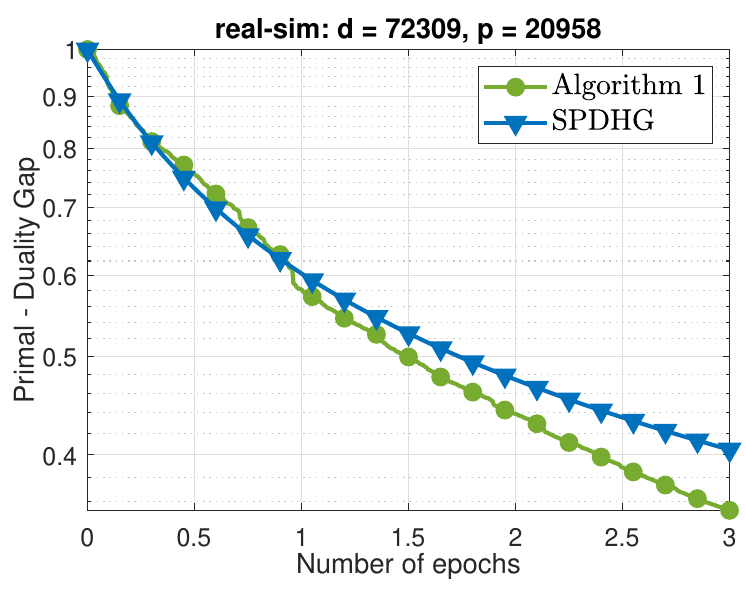}
\hspace{-2ex}
\vspace{-1ex}
\caption{The performance of Algorithm~\ref{alg:A1_main} and SPDHG with single coordinate, i.e., $p_i = 1$ ($i\in [n]$).}
\label{fig:single_coordinates}
\end{center}
\vspace{-2ex}
\end{figure*}

\vspace{1ex}
\noindent\textbf{Experiment 3: Single coordinate update.}
We provide an experiment to test Algorithm \ref{alg:A1_main} and SPDHG using single coordinate (i.e., $p_i = 1$ for all $i \in [n]$, each block has a single entry).
Figure \ref{fig:single_coordinates} shows the performance of two algorithms on the \textbf{w8a}, \textbf{rcv1}, and \textbf{real-sim}  datasets. 
We choose $\rho_0 := 10/\norms{K}$ in Algorithm~\ref{alg:A1_main} and $\tau = \sigma := 10/\norms{K}$ in SPDHG among all datasets.
Since the per-iteration complexity of these algorithms is at most $\BigO{\max\set{p,d}}$, we run these algorithms up to $3p$ and $3m$ iterations, respectively, corresponding to $3$ epochs.
From Figure~\ref{fig:single_coordinates}, we can see that SPDHG  performs better than Algorithm~\ref{alg:A1_main} on the \textbf{w8a} and \textbf{rcv1} datasets. 
However, Algorithm~\ref{alg:A1_main} is better than SPDHG on the \textbf{real-sim} dataset. 

\begin{figure*}[ht!]
\begin{center}
\includegraphics[width=0.48\textwidth]{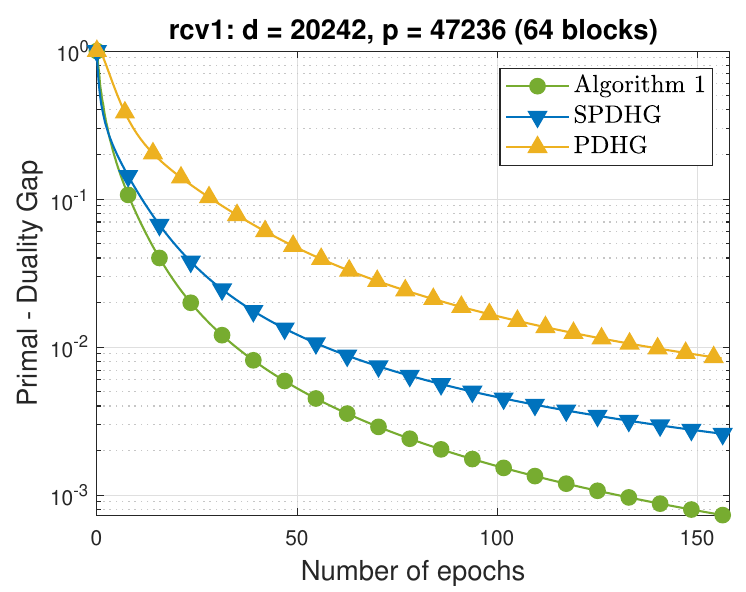}
\includegraphics[width=0.48\textwidth]{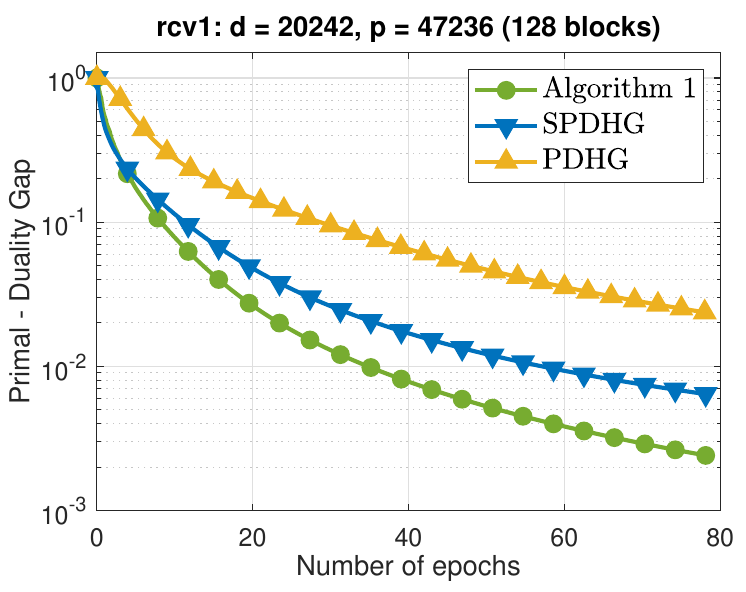}
\includegraphics[width=0.48\textwidth]{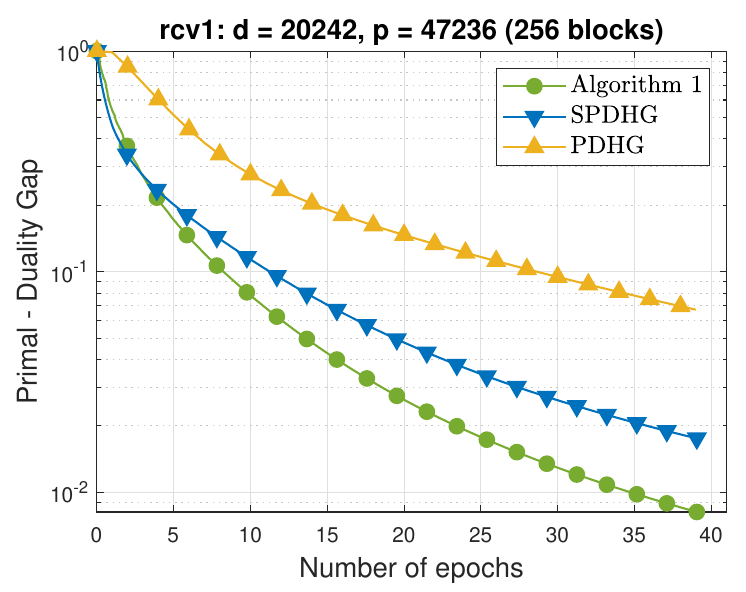}
\includegraphics[width=0.48\textwidth]{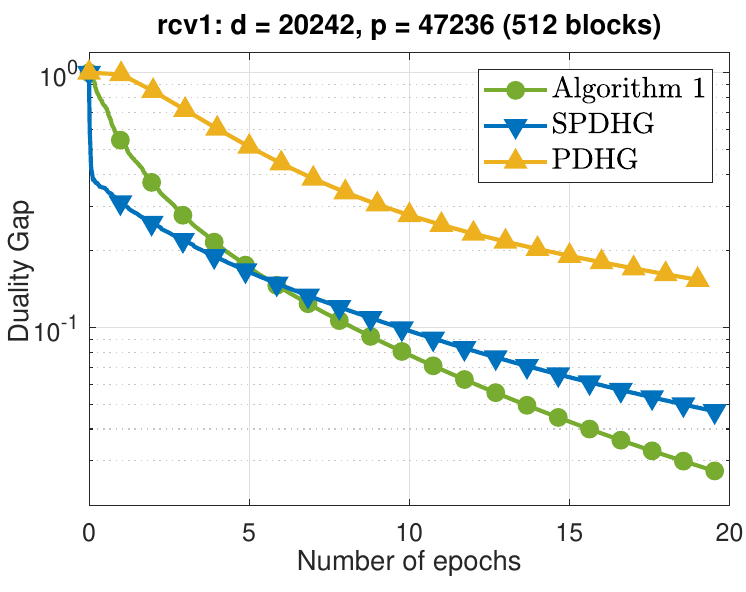}
\vspace{-1ex}
\caption{Comparing Algorithm \ref{alg:A1_main} and SPDHG using different number of blocks: $n = 64$, $128$, $256$, and $512$ on the \textbf{rcv1} dataset.}
\label{fig:different_blocks}
\end{center}
\vspace{-3ex}
\end{figure*}

\vspace{1ex}
\noindent\textbf{Experiment 4: Using different block-coordinate sizes.}
In this experiment, we test the effect of the number of block-coordinates on the performance of Algorithm~\ref{alg:A1_main} and SPDHG.
We still compare them with PDHG. 
We only choose the \textbf{rcv1} dataset since it has relatively large $p$ and $d$ ($d = 20242$ and $p = 47236$).
We choose the number of blocks $n$ to be $64$, $128$, $256$, and $512$.
We choose $\rho_0 := 10/\norms{K}$ in Algorithm~\ref{alg:A1_main},  $\tau = \sigma := 10/\norms{K}$ in SPDHG, and $\tau := 10/\norms{K}, \sigma := 0.03$ in PDHG for all cases.
The performance of three algorithms is shown in Figure~\ref{fig:different_blocks} for a fixed number of iterations.

From Figure~\ref{fig:different_blocks}, we can see that  Algorithm~\ref{alg:A1_main} still performs well and better than SPDHG as well as PDHG.
Hence, Algorithm~\ref{alg:A1_main} seems to work well on \eqref{eq:svm_problem} when running it with block coordinates.

\beforesubsec
\subsection{Least absolute deviation $($LAD$)$ problem}\label{subsec_lad}
\aftersubsec
We consider the following well-studied least absolute deviations (LAD) problem:
\myeq{eq:LAD_exam}{
\min_{x \in \R^p}\Big\{ F(x) := \norms{Kx - b}_1 + \lambda\norms{x}_1 \Big\},
}
where $K\in\R^{d\times p}$, $b\in\R^d$ and $\lambda > 0$ is a regularization parameter.

\begin{figure*}[hpt!]
\begin{center}
\hspace{-2ex}
\includegraphics[width=0.33\textwidth]{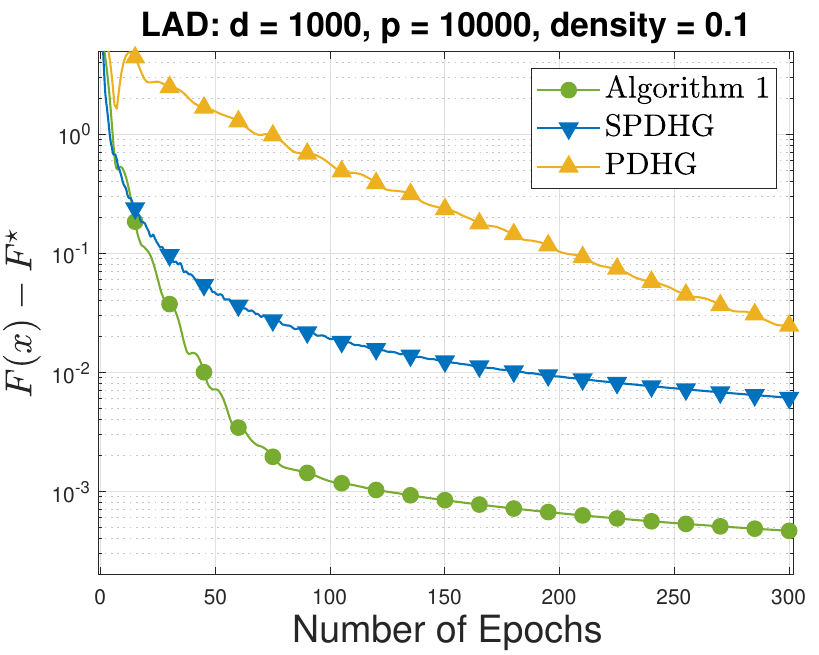}
\includegraphics[width=0.33\textwidth]{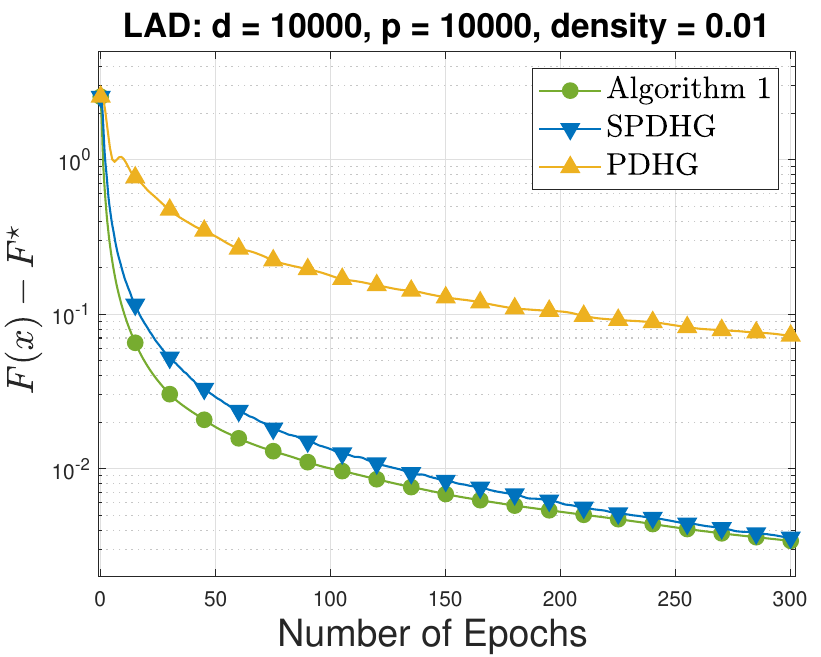}
\includegraphics[width=0.33\textwidth]{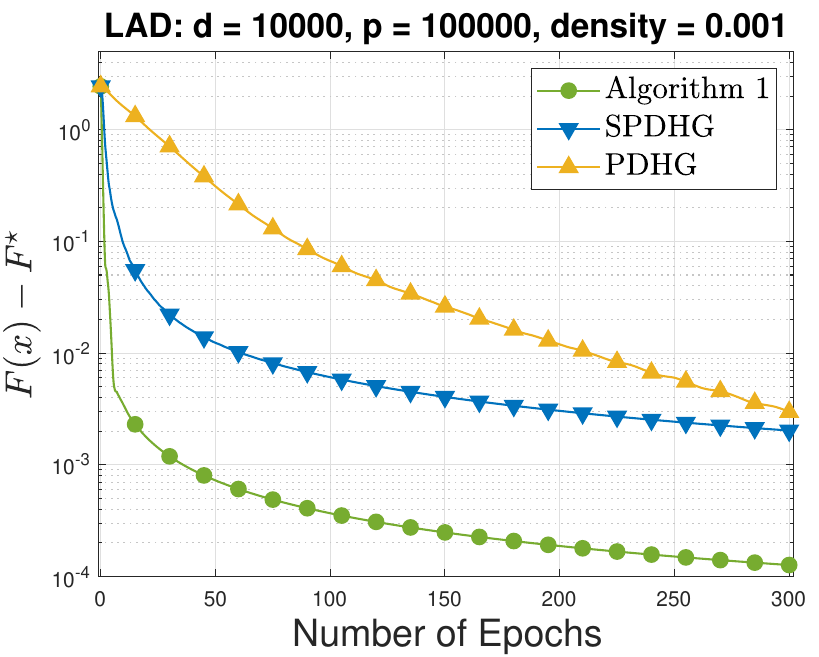}
\hspace{-2ex}
\caption{Comparison of Algorithm \ref{alg:A1_main} with PDHG and SPDHG on \eqref{eq:LAD_exam} using synthetic data.}
\label{fig:LAD}
\end{center}
\vspace{-3ex}
\end{figure*}

We again test Algorithm~\ref{alg:A1_main} and compare it with SPDHG and PDHG on three problem instances, where $K$ is generated from the standard Gaussian distribution with different densities. Here, we choose $\lambda := 1/d$ ($d$ is the number of rows of $K$) and $b := Kx^{\natural} + 0.1\Lc(0,1)$, where $x^{\natural}$ is a predefined sparse vector and $\Lc$ stands for Laplace noise.
The experiment results are reported in Figure~\ref{fig:LAD}, where we run for $300$ epochs and use $32$ blocks in the stochastic algorithms. 

For \eqref{eq:LAD_exam}, we choose $3$ instances, where one case is dense with $10\%$ nonzero entries in $K$, and two other instances are sparse with only $1\%$ and $0.1\%$ nonzero entries, respectively.
After a fining tune, we choose the parameter $\rho_0$ of Algorithm~\ref{alg:A1_main} and the step-size $\tau$ and $\sigma$ for SPDHG and PDHG, respectively as follows.
\begin{itemize}
\item Instance 1 with $10\%$ nonzero entries, we choose $\rho_0 := 10/\norms{K}$ in Algorithm~\ref{alg:A1_main},  $\tau := 0.005, \sigma := 0.01$ in SPDHG, and $\tau := 0.005, \sigma := 0.01$ in PDHG.
\item Instance 2 with $1\%$ nonzero entries, we choose $\rho_0 := 50/\norms{K}$ in Algorithm~\ref{alg:A1_main}),  $\tau := 0.03, \sigma := 0.01$ in SPDHG, and $\tau := 0.005, \sigma := 0.1$ in PDHG.
\item Instance 3 with $0.1\%$ nonzero entries, we choose $\rho_0 := 100/\norms{K}$ in Algorithm~\ref{alg:A1_main},  $\tau := 0.01, \sigma := 0.05$ in SPDHG, and $\tau := 0.001, \sigma := 0.5$ in PDHG.
\end{itemize}
Note that the choice of $\rho_0$ simply trades off the effect of the primal and dual initial points to the complexity bounds as we can see in the right-hand side bounds of Theorem~\ref{th:convergence_rate} and Theorem~\ref{th:convergence_rate2}.

We can observe from Figure~\ref{fig:LAD} that Algorithm~\ref{alg:A1_main} still works well compared to SPDHG under $3$ different instances.
As expected, both Algorithm~\ref{alg:A1_main} and SPDHG outperform PDHG in all cases.

\beforesec
\section{The Proof of The Main Results: Theorem~\ref{th:convergence_rate} and Theorem~\ref{th:convergence_rate2}}\label{apdx:sec:appendix1}
\aftersec
This section provides the full proof of Theorem~\ref{th:convergence_rate} and Theorem~\ref{th:convergence_rate2}.

\beforesubsec
\subsection{Preliminary results}\label{apdx:subsec:proof_tech1}
\aftersubsec
The following identities will be repeatedly used for our convergence analysis.
\begin{enumerate}
\item[$\mathrm{(i)}$] For any $a, b, u\in\R^p$ and $\tau \in [0, 1]$, we have
\myeq{eq:useful_id1}{
\tau(1 - \tau)\norms{u - a}^2 + \norms{(1-\tau)a + \tau u - b}^2  = \tau \norms{u - b}^2 + (1 - \tau)\norms{b - a}^2.
}
\item[$\mathrm{(ii)}$]
 For any $a, \hat{a}\in\R^p$, $\tau\in [0, 1]$, $\rho > 0$, and $\hat{\rho} > 0$, we have
\myeq{eq:useful_id2}{
\hspace{-2ex}
\arraycolsep=0.2em
\begin{array}{ll}
(1 -  \tau)\rho\norms{a  -   \hat{a}}^2 & + {~} \tau\rho\norms{\hat{a}}^2 - (1-\tau)(\rho-\hat{\rho})\norms{a}^2  =   \rho\norms{\hat{a} - (1-\tau)a}^2 \vspace{1ex}\\
& + {~} (1-\tau)[\hat{\rho} - (1-\tau)\rho]\norms{a}^2.
\end{array}
\hspace{-1ex}
}
\item[$\mathrm{(iii)}$]
For any $a, b\in\R^p$, $\rho > 0$, and $\hat{\rho} > \rho$, we have
\myeq{eq:useful_id3}{
\rho\norm{a}^2 - \hat{\rho}\norms{b}^2 \leq \frac{\rho\hat{\rho}}{\hat{\rho} - \rho}\norms{a - b}^2.
}
\end{enumerate}
In addition, the following lemma will also be used in the sequel, whose proof is in \cite{fercoq2015accelerated}.

\begin{lemma}\cite{fercoq2015accelerated}\label{le:vector_x_r_rep}
Given a sequence $\sets{\tilde{x}^k}_{k\geq 0}$ in $\R^p$ and a nonincreasing sequence $\sets{\tau_k}_{k\geq 0}$ in $(0, 1]$, let $\sets{(x^k, \hat{x}^k)}_{k\geq 0}$ be updated as
\myeqn{
\hat{x}^{k} :=   (1-\tau_k)x^k + \tau_k\tilde{x}^k \quad \text{and} \quad x^{k+1} := \hat{x}^k + \tfrac{\tau_k}{\tau_0}(\tilde{x}^{k+1} - \tilde{x}^k).
}
Then, we have
\myeq{eq:x_r_breve}{
\hspace{-2ex}
x^k = \sum_{l=0}^k\gamma_{k,l}\tilde{x}^l, \quad \text{where}\quad
\gamma_{k+1,l} := \begin{cases}
(1-\tau_k)\gamma_{k,l} &\text{if}~l=0,\cdots, k-1\\
(1-\tau_k)\gamma_{k,k} + \tau_k - \frac{\tau_k}{\tau_0} &\text{if}~l=k,\\
\frac{\tau_k}{\tau_0} &\text{if}~l=k+1,
\end{cases}
\hspace{-2ex}
}
and $\gamma_{0,0} := 1$.
Moreover, we have  $\gamma_{k,l} \geq 0$ for $l=0,\cdots, k$ and $\sum_{l=0}^k\gamma_{k,l} = 1$ for $k\geq 0$.
\end{lemma}

\beforesubsec
\subsection{Properties of the augmented Lagrangian function}\label{apdx:subsec:aug_lag_properties}
\aftersubsec
The augmented Lagrangian function $\Lc_{\rho}$ defined by \eqref{eq:aug_Lag} will serve as a \mytbi{merit function} for our convergence analysis.
To investigate the properties of $\Lc_{\rho}(\cdot,\cdot)$, we consider
\myeq{eq:appdix_phi_func}{
\psi_{\rho}(x, w, y) := \iprods{y, Kx + Bw - b} + \frac{\rho}{2}\norms{Kx + Bw - b}^2.
}
This function has following properties:
\myeq{eq:psi_deriv}{
\hspace{-2ex}
\arraycolsep=0.1em
\left\{\begin{array}{ll}
& \nabla_{w}\psi_{\rho}(x, w, y) = B^{\top}(y + \rho(Kx + Bw - b)), \vspace{1ex}\\
& \nabla_{x_i}\psi_{\rho}(x, w, y) = K_{i}^{\top}(y + \rho(Kx + Bw - b)), \vspace{1ex}\\
&\Vert \nabla_{x_i}\psi_{\rho}(x + U_id_i,  w, y) - \nabla_{x_i}\psi_{\rho}(x,  w, y) \Vert =  \rho\Vert K_{i}^{\top}K_{i}d_i\Vert \leq \rho \Vert K_{i}\Vert^2\norms{d_i},~~\forall d_i\in\R^{p_i}.
\end{array}\right.
\hspace{-5ex}
}
These estimates allow us to conclude that $\nabla_{x_i}\psi_{\rho}(x + U_i(\cdot), w, y)$ is Lipschitz continuous with the Lipschitz constant $\rho\norms{K_i}^2$.
Directly using the definition of $\psi_{\rho}$, for all $x$, $w$, $y$, $\hat{x}$, and $\hat{w}$, we also have the following identity:
\myeq{eq:lower_bound}{
\arraycolsep=0.2em
\begin{array}{lcl}
\psi_{\rho}(\hat{x}, \hat{w}, y) & = & \psi_{\rho}(x, w, y) +  \iprods{\nabla_x\psi_{\rho}(x, w, y), \hat{x} - x} + \iprods{\nabla_w\psi_{\rho}(x, w, y), \hat{w} - w} \vspace{1ex}\\
&& + {~} \frac{\rho}{2}\norms{K(\hat{x} - x) + B(\hat{w} - w)}^2.
\end{array}
}
As a consequence of \eqref{eq:lower_bound}, with $\bar{L}_{\sigma}$ defined by \eqref{eq:common_quantities}, if $\hat{x}$ only changes one block $i$ from $x$ to $\hat{x} = x + U_id_i$ for any $i\in [n]$, then  we have
\myeq{eq:lower_bound2}{
\arraycolsep=0.2em
\begin{array}{lcl}
\psi_{\rho}(\hat{x}, w, y) & \leq & \psi_{\rho}(x, w, y) +  \iprods{\nabla_x\psi_{\rho}(x, w, y), \hat{x} - x} + \frac{\rho \bar{L}_{\sigma}}{2}\norms{\hat{x} - x}_{\sigma}^2.
\end{array}
}
The expressions \eqref{eq:lower_bound} and  \eqref{eq:lower_bound2} are key to our analysis in the sequel. 

\beforesubsec
\subsection{Lyapunov function and key estimates}\label{apdx:subsec:le:key_estimate01}
\aftersubsec
%
\noindent\textbf{Lyapunov function.}
Give a sequence of nonnegative real numbers $\sets{\gamma_{k,l}}_{k,l\geq 0}$ and a sequence $\sets{\tilde{x}^k}_{k\geq 0}$ in $\R^p$, let us introduce the following quantities:
\myeq{eq:non_smooth}{
\bar{f}^k_i :=  \displaystyle\sum_{l=0}^k\gamma_{k,l}f_i(\tilde{x}_i^l) \quad \text{and} \quad \bar{f}^k :=  \displaystyle\sum_{i=1}^n\bar{f}_i^k.
}
For $\psi_{\rho}$ defined by \eqref{eq:appdix_phi_func} and $\bar{f}^k$ defined by \eqref{eq:non_smooth}, we also introduce
\myeq{eq:L_hat_k}{
\bar{\Lc}_k(y) := \bar{f}^k  + h(x^k)  + g(w^k) + \psi_{\rho_{k-1}}(x^k, w^k, y).
}
Given \eqref{eq:L_hat_k} and $\Lc$ defined by \eqref{eq:minimax_form}, we define a Lyapunov function as follows:
\myeq{eq:lyapunov_func}{
\arraycolsep=0.2em
\begin{array}{lcl}
\Ec_k(x, w, y) &:= &  \bar{\Lc}_k(y) -  \Lc(x, w, \bar{y}^k) + \frac{1}{2\eta_{k-1}}\norms{\hat{y}^{k} - y}^2   \vspace{1ex}\\
&& + {~} \sum_{i=1}^n\frac{\tau_{k-1}}{2q_i}\left(\frac{\tau_{k-1}\sigma_i}{\tau_0\beta_{k-1}} + \mu_{f_i}\right)\norms{\tilde{x}_i^{k} - x_i}^2.
\end{array}
}
\noindent\textbf{Full update vs. block coordinate update.}
For our convergence analysis, we consider the following full update of $\tilde{x}^k := (\tilde{x}^k_1, \cdots, \tilde{x}^k_n)$:
\myeq{eq:r_full_steps}{
\hspace{-3ex}
\arraycolsep=0.1em
\begin{array}{lcll}
\bar{\tilde{x}}^{k+1}_i & := & {\displaystyle\argmin_{x_i}}\Big\{ f_i(x_i) \! + \! \iprods{\nabla_{x_i}{h}(\hat{x}^k) + \nabla_{x_i}{\psi}_{\rho_k}(\hat{x}^k, w^{k+1}, \hat{y}^k), x_i -  \hat{x}^k_i} \vspace{1ex}\\
&& \hspace{8ex} + {~} \frac{\tau_k\sigma_i}{2\tau_0\beta_k}\norms{x_i - \tilde{x}^k_i}^2 \Big\}, \quad \forall i \in [n].
\end{array}
\hspace{-2.5ex}
}
Then, from \eqref{eq:r_full_steps}, Step~\ref{step:A1-i4} of Algorithm~\ref{alg:A1_main} can be shortly rewritten as
\myeq{eq:coordinate_update}{
\tilde{x}^{k+1}_i  = \begin{cases}
\bar{\tilde{x}}^{k+1}_i &\text{if}~i = i_k, \\
\tilde{x}^k_i &\text{otherwise}.
\end{cases}
}
\beforesubsec
\subsection{Preparatory: Two intermediate steps of convergence analysis}
\aftersubsec
The following two lemmas serve as key estimates for our convergence analysis of Algorithm~\ref{alg:A1_main} in the sequel.

\begin{lemma}\label{le:key_estimate_f_g}
Let $\set{(x^k, \tilde{x}^k, w^{k+1}, \hat{y}^k)}$ be generated by Algorithm~\ref{alg:A1_main}, $\bar{f}^k$ be defined by \eqref{eq:non_smooth}, and $\psi_{\rho}$ be defined by \eqref{eq:appdix_phi_func}.
Then, for any fixed $x \in\dom{f}$, it holds that
\myeq{eq:main_bound_f}{
\hspace{-2ex}
\arraycolsep=0.2em
\begin{array}{l}
\Exps{i_k}{  \bar{f}^{k+1} + \sum_{i=1}^{n}\frac{\tau_k}{2q_i}\big(\frac{\tau_k\sigma_i}{\tau_0\beta_k}  + \mu_{f_i}\big)\norms{\tilde{x}_i^{k+1} - x_i }^2} \leq (1-\tau_k)\bar{f}^{k} + \tau_kf(x) \vspace{1.5ex}\\
\hspace{5ex} + {~} \frac{\tau_k}{\tau_0}\sum_{i=1}^nq_i\big\langle \nabla_{x_i}h(\hat{x}^k)  + \nabla_{x_i}{\psi}_{\rho_k}(\hat{x}^k, w^{k+1}, \hat{y}^k), (1 - \frac{\tau_0}{q_i})\tilde{x}_i^k + \frac{\tau_0}{q_i}x_i - \bar{\tilde{x}}^{k+1}_i \big\rangle \vspace{1ex}\\
\hspace{5ex} + {~} \sum_{i=1}^n \frac{\tau_k}{2q_i}\big[ \frac{\tau_k\sigma_i}{\tau_0\beta_k} + (1 - q_i)\mu_{f_i}\big] \norms{\tilde{x}_i^k - x_i}^2 -  \frac{\tau_k^2}{2\tau_0^2\beta_k}\sum_{i=1}^n\sigma_iq_i\norms{\bar{\tilde{x}}^{k+1}_i - \tilde{x}_i^k}^2, \vspace{1.5ex}\\
\Exps{i_k}{\psi_{\rho_k}(x^{k+1}, w^{k+1}, \hat{y}^k)} \leq \psi_{\rho_k}(\hat{x}^k, w^{k+1},\hat{y}^k) + \frac{\rho_k\tau_k^2\bar{L}_{\sigma}}{2\tau_0^2}\sum_{i=1}^nq_i\sigma_i\norms{\bar{\tilde{x}}^{k+1}_i - \tilde{x}_i^k}^2 \vspace{1ex}\\
\hspace{5ex} + {~}  \frac{\tau_k}{\tau_0}\sum_{i=1}^nq_i\iprods{\nabla_{x_i}\psi_{\rho_k}(\hat{x}^k, w^{k+1}, \hat{y}^k), \bar{\tilde{x}}^{k+1}_i - \tilde{x}_i^k}, \vspace{1.5ex}\\
\Exps{i_k}{h(x^{k+1})}  \leq  h(\hat{x}^k) + \sum_{i=1}^n  \frac{\tau_kq_i}{\tau_0} \Big[ \iprods{\nabla_{x_i}{h}(\hat{x}^k), \bar{\tilde{x}}^{k+1}_i - \tilde{x}_i^k} + \frac{\tau_kL_{\sigma}^h\sigma_i}{2\tau_0}\norms{\bar{\tilde{x}}_i^{k+1} - \tilde{x}_i^k}^2 \Big].
\end{array}
\hspace{-4ex}
}
\end{lemma}

\begin{proof}
First, the optimality condition of \eqref{eq:r_full_steps} for $x$ can be read as
\myeq{eq:opt_cond_x}{
\begin{array}{ll}
0 = \nabla{f_i}(\bar{\tilde{x}}^{k+1}_i) +  \nabla_{x_i}h(\hat{x}^k) +  \nabla_{x_i}{\psi}_{\rho_k}(\hat{x}^k, w^{k+1},\hat{y}^k) + \frac{\tau_k\sigma_i}{\tau_0\beta_k}(\bar{\tilde{x}}^{k+1}_i - \tilde{x}^k_i),
\end{array}
}
for some $\nabla{f_i}(\bar{\tilde{x}}^{k+1}_i) \in \partial{f_i}(\bar{\tilde{x}}^{k+1}_i)$.

By $\mu_{f_i}$-convexity of $f_i$, \eqref{eq:opt_cond_x}, for any $\breve{x}_i\in\R^{p_i}$, we can derive
\myeq{eq:f_i_1}{
\begin{array}{lcl}
f_i(\bar{\tilde{x}}^{k+1}_i) &\leq & f_i(\breve{x}_i) + \iprods{\nabla{f_i}(\bar{\tilde{x}}^{k+1}_i), \bar{\tilde{x}}^{k+1}_i - \breve{x}_i} - \frac{\mu_{f_i}}{2}\norms{\bar{\tilde{x}}^{k+1}_i - \breve{x}_i}^2  \vspace{1ex}\\
& \overset{\tiny \eqref{eq:opt_cond_x}}{=} & f_i(\breve{x}_i) + \iprods{ \nabla_{x_i}h(\hat{x}^k)  + \nabla_{x_i}{\psi}_{\rho_k}(\hat{x}^k, w^{k+1},\hat{y}^k), \breve{x}_i - \bar{\tilde{x}}^{k+1}_i} \vspace{1ex}\\
&& + {~} \frac{\tau_k\sigma_i}{\tau_0\beta_k}\iprods{\bar{\tilde{x}}^{k+1}_i - \tilde{x}^k_i, \breve{x}_i - \bar{\tilde{x}}^{k+1}_i} - \frac{\mu_{f_i}}{2}\norms{\bar{\tilde{x}}^{k+1}_i - \breve{x}_i}^2.
\end{array}
}
Now, using $\breve{x}_i := (1 - \frac{\tau_0}{q_i})\tilde{x}_i^k + \frac{\tau_0}{q_i}x_i$ with $\frac{\tau_0}{q_i} \in [0, 1]$ and  $2\iprods{a,b} = \norms{a+b}^2 - \norms{a}^2 - \norms{b}^2$, it is easy to show that
\myeq{eq:key_estimate_f_g_est1}{
\hspace{-2ex}
\arraycolsep=0.2em
\begin{array}{lcl}
\iprods{\bar{\tilde{x}}^{k+1}_i - \tilde{x}^k_i, \breve{x}_i - \bar{\tilde{x}}^{k+1}_i} & = & \iprods{\bar{\tilde{x}}^{k+1}_i - \tilde{x}^k_i, (1 - \frac{\tau_0}{q_i})(\tilde{x}_i^k - \bar{\tilde{x}}^{k+1}_i) + \frac{\tau_0}{q_i}(x_i - \bar{\tilde{x}}^{k+1}_i)} \vspace{1ex}\\
&= & \frac{\tau_0}{q_i}\iprods{\bar{\tilde{x}}^{k+1}_i - \tilde{x}^k_i, x_i - \bar{\tilde{x}}^{k+1}_i} - (1 - \frac{\tau_0}{q_i})\norms{\tilde{x}_i^k - \bar{\tilde{x}}^{k+1}_i}^2 \vspace{1ex}\\
&\leq & \frac{\tau_0}{2q_i}\norms{x_i - \tilde{x}^k_i}^2  - \frac{\tau_0}{2q_i}\norms{x_i - \bar{\tilde{x}}^{k+1}_i}^2 - \frac{1}{2}\norms{\tilde{x}_i^k - \bar{\tilde{x}}^{k+1}_i}^2.
\end{array}
\hspace{-2ex}
}
Again, by $\mu_{f_i}$-convexity of $f_i$, we can deduce that
\myeq{eq:key_estimate_f_g_est2}{
\hspace{-2ex}
\arraycolsep=0.1em
\hspace{-0ex}\begin{array}{lcl}
f_i(\breve{x}_i) - \frac{\mu_{f_i}}{2}\norms{\bar{\tilde{x}}_i^{k+1} - \breve{x}_i}^2  & \leq &  \big(1 - \frac{\tau_0}{q_i}\big)f_i(\tilde{x}_i^k) + \frac{\tau_0}{q_i}f_i(x_i)  - \frac{\mu_{f_i}}{2} \big(1 - \frac{\tau_0}{q_i}\big) \frac{\tau_0}{q_i}\norms{x_i - \tilde{x}_i^k}^2 \vspace{1ex}\\
&& - {~} \frac{\mu_{f_i}}{2}\norms{\big(1 - \frac{\tau_0}{q_i}\big)\tilde{x}_i^k + \frac{\tau_0}{q_i}x_i - \bar{\tilde{x}}_i^{k+1}}^2 \vspace{1ex}\\
&  \overset{\tiny\eqref{eq:useful_id1}}{=} &  \big(1 - \frac{\tau_0}{q_i}\big)f_i(\tilde{x}_i^k) + \frac{\tau_0}{q_i}f_i(x_i) \vspace{1ex}\\
&& - {~} \frac{\mu_{f_i}}{2}\Big[\frac{\tau_0}{q_i}\norms{\bar{\tilde{x}}_i^{k+1} - x_i}^2 + \big(1 - \frac{\tau_0}{q_i}\big)\norms{\bar{\tilde{x}}_i^{k+1} - \tilde{x}^k_i}^2 \Big] \vspace{1ex}\\
& \leq  &  \big(1 - \frac{\tau_0}{q_i}\big)f_i(\tilde{x}_i^k) + \frac{\tau_0}{q_i}f_i(x_i) - \frac{\tau_0\mu_{f_i}}{2q_i}\norms{\bar{\tilde{x}}_i^{k+1} - x_i}^2. 
\end{array}
\hspace{-4ex}
}
Therefore, plugging \eqref{eq:key_estimate_f_g_est1} and \eqref{eq:key_estimate_f_g_est2} into \eqref{eq:f_i_1}, and  using again $\breve{x}_i := (1 - \frac{\tau_0}{q_i})\tilde{x}_i^k + \frac{\tau_0}{q_i}x_i$, we can further derive
\myeq{eq:key_estimate_f_g_est3}{
\arraycolsep=0.2em
\begin{array}{lcl}
f_i(\bar{\tilde{x}}^{k+1}_i) & \leq & \big(1 - \frac{\tau_0}{q_i}\big)f_i(\tilde{x}_i^k) + \frac{\tau_0}{q_i}f_i(x_i) - \frac{\tau_0\mu_{f_i}}{2q_i}\norms{\bar{\tilde{x}}_i^{k+1} - x_i}^2  \vspace{1ex}\\
&& + {~} \frac{\tau_k\sigma_i}{2q_i\beta_k}\left[\norms{x_i - \tilde{x}^k_i}^2 - \norms{x_i - \bar{\tilde{x}}_i^{k+1}}^2\right] - \frac{\tau_k\sigma_i}{2\tau_0\beta_k}\norms{\tilde{x}_i^k - \bar{\tilde{x}}^{k+1}_i}^2 \vspace{1.5ex}\\
&& + {~} \big\langle{ \nabla_{x_i}h(\hat{x}^k)  + \nabla_{x_i}{\psi}_{\rho_k}(\hat{x}^k, w^{k+1},\hat{y}^k), (1 - \frac{\tau_0}{q_i})\tilde{x}_i^k + \frac{\tau_0}{q_i}x_i - \bar{\tilde{x}}^{k+1}_i}\big\rangle.
\end{array}
}
Next, using \eqref{eq:x_r_breve} of Lemma~\ref{le:vector_x_r_rep} into \eqref{eq:non_smooth}, we can show that
\myeqn{
\begin{array}{ll}
\bar{f}^{k+1} &:= \sum_{l=0}^{k+1}\gamma_{k+1,l}f(\tilde{x}^l) \overset{\tiny\eqref{eq:x_r_breve}}{=} (1-\tau_k)\bar{f}^k + \tau_kf(\tilde{x}^k) + \frac{\tau_k}{\tau_0}\left[ f(\tilde{x}^{k+1}) - f(\tilde{x}^k)\right].
\end{array}
}
Taking conditional expectation $\Exps{i_k}{\cdot}$ of this expression, we can further derive
\myeqn{
\arraycolsep=-0.1em
\begin{array}{lcl}
\Exps{i_k}{\bar{f}^{k+1} } & = &  (1-\tau_k)\bar{f}^k + \tau_kf(\tilde{x}^k) + \frac{\tau_k}{\tau_0}  \sum_{i=1}^nq_i\left[ f_i(\bar{\tilde{x}}_i^{k+1})  - f_i(\tilde{x}^k_i)\right] \vspace{1ex}\\
& \overset{\eqref{eq:key_estimate_f_g_est3}}{\leq} & (1-\tau_k)\bar{f}^k + \tau_kf(x) + \frac{\tau_k^2}{2\tau_0\beta_k}\sum_{i=1}^n\sigma_i\left[\norms{x_i - \tilde{x}^k_i}^2 - \norms{x_i - \bar{\tilde{x}}_i^{k+1}}^2\right]\vspace{1ex}\\
&& - {~} \frac{\tau_k}{2}\sum_{i=1}^n\mu_{f_i}\norms{\bar{\tilde{x}}_i^{k+1} - x_i}^2   - \frac{\tau_k^2}{2\tau_0^2\beta_k}\sum_{i=1}^n\sigma_iq_i\norms{\tilde{x}_i^k - \bar{\tilde{x}}^{k+1}_i}^2 \vspace{1ex}\\
&& + {~} \frac{\tau_k}{\tau_0}\sum_{i=1}^nq_i\big\langle{ \nabla_{x_i}h(\hat{x}^k)  + \nabla_{x_i}{\psi}_{\rho_k}(\hat{x}^k, w^{k+1},\hat{y}^k), (1 - \frac{\tau_0}{q_i})\tilde{x}_i^k + \frac{\tau_0}{q_i}x_i - \bar{\tilde{x}}^{k+1}_i}\big\rangle.
\end{array}
}
Finally, substituting the following expressions
\myeqn{
\arraycolsep=0.1em
\left\{\begin{array}{ll}
&\Exps{i_k}{\sum_{i=1}^n\frac{\sigma_i}{q_i}\big( \norms{\tilde{x}^k_i -  x_i}^2 - \norms{\tilde{x}_i^{k+1} -  x_i}^2 \big) }  =  \sum_{i=1}^n\sigma_i[\norms{\tilde{x}^k_i - x_i}^2 - {~} \norms{\bar{\tilde{x}}_i^{k+1} -  x_i}^2], \vspace{1ex}\\
&\Exps{i_k}{\sum_{i=1}^n\frac{\mu_{f_i}}{q_i} \big( \norms{\tilde{x}_i^{k+1} -  x_i}^2 - (1-q_i)\norms{\tilde{x}_i^k - x_i}^2 \big)} = \sum_{i=1}^n\mu_{f_i}\norms{\bar{\tilde{x}}_i^{k+1} -  x_i}^2,
\end{array}\right.
}
into the last inequality and rearranging the result we eventually obtain the first estimate  of \eqref{eq:main_bound_f}. 

Since $x^{k+1} := \hat{x}^k + \frac{\tau_k}{\tau_0}(\tilde{x}^{k+1} - \tilde{x}^k)$, it is clear that $x^{k+1}$ is different than $\hat{x}^k$ only in one block $i_k$. 
By utilizing \eqref{eq:lower_bound2}, we obtain
\myeq{eq:main_estimate_psi_h_est1}{
\arraycolsep=0.2em
\begin{array}{lcl}
\psi_{\rho_k}(x^{k+1}, w^{k+1}, \hat{y}^k) & \leq  & \psi_{\rho_k}(\hat{x}^k, w^{k+1}, \hat{y}^k)  + \iprods{\nabla_{x}\psi_{\rho_k}(\hat{x}^k, w^{k+1}, \hat{y}^k), x^{k+1} - \hat{x}^k} \vspace{1ex}\\
&& + {~} \frac{\rho_k\bar{L}_{\sigma}}{2}\sum_{i=1}^n\sigma_i\norms{x_i^{k+1} - \hat{x}_i^k}^2.
\end{array}
}
Alternatively, by \eqref{eq:Lh_smooth} and using $L_{\sigma}^h$ defined by \eqref{eq:common_quantities}, we also have
\myeq{eq:main_estimate_psi_h_est2}{
h(x^{k+1}) \leq h(\hat{x}^k) + \iprods{\nabla_{x}{h}(\hat{x}^k), x^{k+1} - \hat{x}^k} + \frac{L_{\sigma}^h}{2}\sum_{i=1}^n\sigma_i\norms{x_i^{k+1} - \hat{x}_i^k}^2.
}
Next, by Step~\ref{step:A1-i4} of Algorithm~\ref{alg:A1_main} and \eqref{eq:coordinate_update}, one can establish that
\myeqn{
\arraycolsep=0.2em
\left\{\begin{array}{lcl}
\Exps{i_k}{\iprods{\nabla_{x}\psi_{\rho_k}(\hat{x}^k, w^{k+1}, \hat{y}^k), x^{k+1} - \hat{x}^k} } &= & \frac{\tau_k}{\tau_0}\sum_{i=1}^nq_i\iprods{\nabla_{x_i}\psi_{\rho_k}(\hat{x}^k, w^{k+1}, \hat{y}^k),\bar{\tilde{x}}^{k+1}_i - \tilde{x}^k_i}, \vspace{1.25ex}\\
\Exps{i_k}{\iprods{\nabla_{x}{h}(\hat{x}^k), x^{k+1} - \hat{x}^k} } &= &  \frac{\tau_k}{\tau_0}\sum_{i=1}^nq_i\iprods{\nabla_{x_i}{h}(\hat{x}^k), \bar{\tilde{x}}^{k+1}_i - \tilde{x}_i^k}.
\end{array}\right.
}
Finally, taking conditional expectation of \eqref{eq:main_estimate_psi_h_est1} and \eqref{eq:main_estimate_psi_h_est2}, and substituting the above equalities into the results, 
we obtain the last two estimates of \eqref{eq:main_bound_f}. 
\end{proof}

\begin{lemma}\label{le:key_estimate01}
Let $\set{(x^k, \tilde{x}^k, w^k, \hat{y}^k)}$ be generated by Algorithm~\ref{alg:A1_main} and $\bar{\Lc}_k(\cdot)$ be defined by \eqref{eq:L_hat_k}.
Then, for any $(x, w)\in\dom{F}$, the following estimate holds:
\myeq{eq:main_bound01}{
\hspace{-0.5ex}
\arraycolsep=0.1em
\begin{array}{ll}
&\Exps{i_k}{\bar{\Lc}_{k+1}(\hat{y}^k)  +  \sum_{i=1}^n \frac{\tau_k}{2q_i}\big(\frac{\tau_k\sigma_i}{\tau_0\beta_k}  +  \mu_{f_i}\big)\norms{\tilde{x}_i^{k+1}  -  x_i}^2 }  \leq (1-\tau_k)\bar{\Lc}_k(\hat{y}^k) \vspace{1ex}\\
&\hspace{15ex} + {~} \tau_k\big[ F(x, w) +  \iprods{\hat{y}^k + \rho_k(K\hat{x}^k + Bw^{k+1} - b), Kx + Bw - b} \big] \vspace{1ex}\\
&\hspace{15ex} + {~} \sum_{i=1}^n \frac{\tau_k}{2q_i}\big[ \frac{\tau_k\sigma_i}{\tau_0\beta_k} + (1-q_i)\mu_{f_i}\big] \norms{\tilde{x}_i^k - x_i}^2  \vspace{1ex}\\
&\hspace{15ex} - {~} \frac{\tau_k^2}{2\tau_0^2}\left( \frac{1}{\beta_k}  - \rho_k\bar{L}_{\sigma} -  L^h_{\sigma} \right)\sum_{i=1}^n\sigma_iq_i\norms{\bar{\tilde{x}}^{k+1}_i - \tilde{x}_i^k}^2 \vspace{1ex}\\
&\hspace{15ex} - {~} \frac{\rho_k}{2} \norms{(K\hat{x}^k + Bw^{k+1} - b) - (1-\tau_k)(Kx^k + Bw^k - b)}^2 \vspace{1ex}\\
&\hspace{15ex} - {~} \frac{(1-\tau_k)}{2}\big[\rho_{k-1} - (1-\tau_k)\rho_k\big] \norms{Kx^k + Bw^k - b}^2.
\end{array}\hspace{-2ex}
}
\end{lemma}

\begin{proof}
First, we write down the optimality condition of \eqref{eq:w_step1} as follows:
\myeqn{
0 \in \partial{g}(w^{k+1}) + B^{\top}( \hat{y}^k + \rho_k (K\hat{x}^k + Bw^{k+1} - b)).
}
Using this condition, the convexity of $g$, and $\nabla_w\psi_{\rho_k}(\hat{x}^k, w^{k+1}, \hat{y}^k) = B^{\top}(\hat{y}^k + \rho_k (K\hat{x}^k + Bw^{k+1} - b))$, for any $\breve{w} := (1-\tau_k)w^k + \tau_kw$ with $w\in\dom{g}$, we have 
\myeqn{
\arraycolsep=0.0em
\begin{array}{lcl}
g(w^{k+1}) & \leq & g(\breve{w}) + \iprods{B^{\top}( \hat{y}^k + \rho_k (K\hat{x}^k + Bw^{k+1} - b)), \breve{w} - w^{k+1}} \vspace{1ex}\\
&\overset{\tiny\eqref{eq:psi_deriv}}{\leq} & (1-\tau_k)g(w^k) +  \tau_kg(w) + \iprods{\nabla_w\psi_{\rho_k}(\hat{x}^k, w^{k+1}, \hat{y}^k), (1-\tau_k)w^k + \tau_kw - w^{k+1}}.
\end{array}
}
Combining the last inequality and \eqref{eq:main_bound_f}, and then using the definition of $\bar{\Lc}_{k}$, we have
\myeq{eq:main_bound01_est1}{
\arraycolsep=0.2em
\hspace{-2ex}\begin{array}{ll}
&\Exps{i_k}{ \bar{\Lc}_{k+1}(\hat{y}^k)  +  \sum_{i=1}^n \frac{\tau_k}{2q_i}\big(\frac{\tau_k\sigma_i}{\tau_0\beta_k}  + \mu_{f_i}\big)\norms{\tilde{x}_i^{k+1}  - x_i}^2 }  \leq  (1-\tau_k)\big[ \bar{f}^k + g(w^k) \big]  \vspace{1ex}\\
&\hspace{15ex}  + {~} \tau_k\big[ f(x) + g(w)\big]  +  \sum_{i=1}^n \frac{\tau_k}{2q_i}\big[  \frac{\tau_k\sigma_i}{\tau_0\beta_k} + (1 - q_i)\mu_{f_i}\big] \norms{\tilde{x}_i^k - x_i}^2 \vspace{1ex}\\
&\hspace{15ex} - {~} \frac{\tau_k^2}{2\tau_0^2}\left( \frac{1}{\beta_k}  - \rho_k\bar{L}_{\sigma} -  L^h_{\sigma} \right)\sum_{i=1}^n\sigma_iq_i\norms{\bar{\tilde{x}}^{k+1}_i - \tilde{x}_i^k}^2\vspace{1ex}\\
&\hspace{15ex} + {~} \psi_{\rho_k}(\hat{x}^k, w^{k+1}, \hat{y}^k) + \tau_k\iprods{ \nabla_{x}{\psi}_{\rho_k}(\hat{x}^k, w^{k+1}, \hat{y}^k), x - \tilde{x}^{k}} \vspace{1ex}\\
&\hspace{15ex} + {~} \iprods{\nabla_w\psi_{\rho_k}(\hat{x}^k, w^{k+1}, \hat{y}^k), (1-\tau_k)w^k + \tau_kw - w^{k+1}} \vspace{1ex}\\
&\hspace{15ex} + {~}  h(\hat{x}^k) + \tau_k\iprods{ \nabla_{x}h(\hat{x}^k)  , x - \tilde{x}^{k}}.
\end{array}\hspace{-6ex}
}
Next, since $\psi_{\rho_k}(x, w, \hat{y}^k) = \iprods{\hat{y}^k, Kx + Bw - b} + \frac{\rho_k}{2}\norms{Kx + Bw - b}^2$, we have
\myeq{eq:main_bound01_est012}{
\hspace{-3ex}
\arraycolsep=0.1em
\begin{array}{lcl}
\Tc_{[1]} & := & \tau_k\psi_{\rho_k}(x, w, \hat{y}^k) -  \frac{\rho_k\tau_k}{2}\norms{(K(x - \hat{x}^k) + B(w - w^{k+1})}^2 \vspace{1ex}\\
&=& \tau_k\iprods{\hat{y}^k + \rho_k(K\hat{x}^k + Bw^{k+1} - b), Kx + Bw - b} - \frac{\rho_k\tau_k}{2}\norms{K\hat{x}^k + Bw^{k+1} - b}^2.  
\end{array}
\hspace{-3ex}
}
Moreover, by Step~\ref{step:A1-i2a} of Algorithm~\ref{alg:A1_main},  we have $\tau_k(x - \tilde{x}^k) = (1-\tau_k)(x^k - \hat{x}^k) + \tau_k(x - \hat{x}^k)$.
Using this expression, we can deduce that
\myeqn{
\hspace{-2ex}
\arraycolsep=0.2em
\begin{array}{lcl}
\Tc_{[2]} &:= & \psi_{\rho_k}(\hat{x}^k, w^{k+1}, \hat{y}^k) +  \tau_k\iprods{\nabla_x\psi_{\rho_k}(\hat{x}^k, w^{k+1}, \hat{y}^k), x - \tilde{x}^k}\vspace{1.25ex}\\
&& + {~} \iprods{\nabla_w\psi_{\rho_k}(\hat{x}^k, w^{k+1}, \hat{y}^k), (1-\tau_k)w^k + \tau_kw - w^{k+1}} \vspace{1ex}\\
& = &  (1-\tau_k)\big[\psi_{\rho_k}(\hat{x}^k, w^{k+1}, \hat{y}^k) + \iprods{\nabla_x\psi_{\rho_k}(\hat{x}^k, w^{k+1}, \hat{y}^k), x^k - \hat{x}^k} \vspace{1ex}\\
&& + {~}  \iprods{\nabla_w\psi_{\rho_k}(\hat{x}^k, w^{k+1}, \hat{y}^k), w^k - w^{k+1}} \big] \vspace{1.25ex}\\
& & + {~} \tau_k\big[\psi_{\rho_k}(\hat{x}^k, w^{k+1}, \hat{y}^k) + \iprods{\nabla_x\psi_{\rho_k}(\hat{x}^k, w^{k+1}, \hat{y}^k), x  - \hat{x}^k} \vspace{1ex}\\
&& + {~}  \iprods{\nabla_w\psi_{\rho_k}(\hat{x}^k, w^{k+1}, \hat{y}^k), w - w^{k+1}} \big].
\end{array}
\hspace{-8ex}
}
Furthermore, utilizing \eqref{eq:main_bound01_est012} and \eqref{eq:lower_bound}, we can further estimate $\Tc_{[2]}$ as
\myeq{eq:main_bound01_est2}{
\hspace{-2ex}
\arraycolsep=0.0em
\begin{array}{lcl}
\Tc_{[2]} 
& \overset{\tiny\eqref{eq:lower_bound}}{=} & (1-\tau_k)\psi_{\rho_k}(x^k, w^{k}, \hat{y}^k) - \frac{(1-\tau_k)\rho_k}{2}\norms{K(x^k - \hat{x}^k) + B(w^k - w^{k+1})}^2\vspace{1.5ex}\\
& & + {~} \tau_k\psi_{\rho_k}(x, w, \hat{y}^k) - \frac{\rho_k\tau_k}{2}\norms{K(x - \hat{x}^k) + B(w - w^{k+1})}^2 \vspace{1ex}\\
& \overset{\tiny\eqref{eq:main_bound01_est012}}{=} & (1-\tau_k)\psi_{\rho_{k-1}}(x^k, w^k, \hat{y}^k) - \frac{(1-\tau_k)\rho_k}{2}\norms{K(x^k - \hat{x}^k) + B(w^k - w^{k+1})}^2  \vspace{1ex}\\
&&  + {~}  \tau_k\iprods{\hat{y}^k + \rho_k(K\hat{x}^k + Bw^{k+1} - b), Kx + Bw - b} - \frac{\rho_k\tau_k}{2}\norms{K\hat{x}^k + Bw^{k+1} - b}^2 \vspace{1.5ex}\\
&&  + {~} \frac{(1-\tau_k)(\rho_k - \rho_{k-1})}{2}\norms{Kx^k + Bw^k - b}^2 \vspace{1ex}\\
& \overset{\eqref{eq:useful_id2}}{=} & (1-\tau_k)\psi_{\rho_{k-1}}(x^k, w^k, \hat{y}^k) + \tau_k\iprods{\hat{y}^k + \rho_k(K\hat{x}^k + Bw^{k+1} - b), Kx + Bw - b} \vspace{1ex}\\
&& - {~} \frac{\rho_k}{2}\norms{K\hat{x}^k + Bw^{k+1} - b - (1-\tau_k)(Kx^k + Bw^k - b)}^2 \vspace{1ex}\\
&& - {~} \frac{(1-\tau_k)}{2}\big[\rho_{k-1} - (1-\tau_k)\rho_k\big]\norms{Kx^k + Bw^{k} - b}^2.
\end{array}
\hspace{-3ex}
}
In addition, we also have
\myeq{eq:main_bound01_est3}{
\hspace{-1ex}
\arraycolsep=0.2em
\begin{array}{lcl}
h(\hat{x}^k) + \tau_k\iprods{\nabla_{x}{h}(\hat{x}^k), x - \tilde{x}^k}  & \leq & h\big((1-\tau_k)x^k + \tau_kx \big)  \leq  (1-\tau_k)h(x^k) + \tau_kh(x).
\end{array}
\hspace{-1ex}
}
Substituting \eqref{eq:main_bound01_est2} and \eqref{eq:main_bound01_est3} into \eqref{eq:main_bound01_est1}, and then simplifying the result, we eventually get \eqref{eq:main_bound01}.
\end{proof}

\beforesubsec
\subsection{Key estimate for Algorithm~\ref{alg:A1_main}}\label{subsec:key_estimate_for_A1}
\aftersubsec
Next, we further estimate \eqref{eq:main_bound01} in terms of $y$ in the following lemma.

\begin{lemma}\label{le:key_spd_scheme3}
Let $\set{(x^k, \tilde{x}^k, w^k, \hat{y}^k, \bar{y}^k)}$ be  generated by Algorithm~\ref{alg:A1_main}, $\Lc$ be defined by \eqref{eq:minimax_form}, and $\bar{\Lc}_k(\cdot)$ be defined by \eqref{eq:L_hat_k}.
Then, for any $(x, w, y) \in \dom{F}\times\R^d$, we have
\myeq{eq:est_1001}{
\hspace{-0.5ex}
\arraycolsep=0.2em
\begin{array}{ll}
& \Exps{i_k}{ \bar{\Lc}_{k+1}(y) - \Lc(x, w, \bar{y}^{k+1}) } \leq (1-\tau_k) \big[ \bar{\Lc}_k(y)   - \Lc(x, w, \bar{y}^k) \big] \vspace{1.5ex}\\
&\hspace{15ex} + {~}  \frac{1}{2\eta_k}\norms{\hat{y}^k - y}^2 - \frac{1}{2\eta_k}\Exps{i_k}{\norms{\hat{y}^{k+1} - y}^2} \vspace{1.5ex}\\
&\hspace{15ex} + {~}  \sum_{i=1}^n\frac{\tau_k}{2q_i}\big[\frac{\tau_k\sigma_i}{\tau_0\beta_k}  + (1-q_i)\mu_{f_i} \big] \norms{\tilde{x}_i^k - x_i}^2\vspace{1ex}\\
&\hspace{15ex} - {~} \Exps{i_k}{\sum_{i=1}^n\frac{\tau_k}{2q_i}\big(\frac{\tau_k\sigma_i}{\tau_0\beta_k}  + \mu_{f_i} \big)\norms{\tilde{x}_i^{k+1} - x_i}^2 }\vspace{1.5ex}\\
&\hspace{15ex} - {~} \frac{\tau_k^2}{2\tau_0^2}\left(\frac{1}{\beta_k} - \rho_k\bar{L}_{\sigma} - L^h_{\sigma} - \frac{\rho_k\eta_k\bar{L}_{\sigma}}{\rho_k-\eta_k}\right)\sum_{i=1}^n\sigma_iq_i\norms{\bar{\tilde{x}}^{k+1}_i - \tilde{x}_i^k}^2\vspace{1.5ex}\\
&\hspace{15ex} - {~} \frac{(1-\tau_k)}{2}\big[\rho_{k-1} - (1-\tau_k)\rho_k\big] \norms{Kx^k + Bw^{k} - b}^2.
\end{array}
\hspace{-5ex}
}
\end{lemma}

\begin{proof}
From \eqref{eq:L_hat_k}, for any $y$, we have $\bar{\Lc}_k(\hat{y}^k) = \bar{\Lc}_k(y) + \iprods{ \hat{y}^k - y, Kx^k + Bw^k - b}$.
Therefore, using the update of $\hat{y}^{k+1}$ from Algorithm~\ref{alg:A1_main}, we can show that
\myeqn{
\arraycolsep=0.1em
\begin{array}{ll}
\bar{\Lc}_{k+1}(\hat{y}^k) & - {~} (1-\tau_k)\bar{\Lc}_k(\hat{y}^k) =  \bar{\Lc}_{k+1}(y) - (1-\tau_k)\bar{\Lc}_k(y)  \vspace{1ex}\\
& + {~} \iprods{\hat{y}^k - y, Kx^{k+1} + Bw^{k+1} - b - (1-\tau_k)(Kx^k + Bw^k - b)} \vspace{1.2ex}\\
&\overset{\tiny \textrm{Step~\ref{step:A1-i6}}}{=}  \bar{\Lc}_{k+1}(y) - (1-\tau_k)\bar{\Lc}_k(y)  + \frac{1}{\eta_k}\iprods{\hat{y}^k - y, \hat{y}^{k+1} - \hat{y}^k} \vspace{1.2ex}\\
& =  \bar{\Lc}_{k+1}(y) - (1-\tau_k)\bar{\Lc}_k(y)  \vspace{1ex}\\
& - {~}  \frac{1}{2\eta_k}\left[\norms{\hat{y}^k - y}^2 - \norms{\hat{y}^{k+1} - y}^2 + \norms{\hat{y}^{k+1} - \hat{y}^k}^2\right].
\end{array}
}
Moreover, since $\bar{y}^{k+1} := (1-\tau_k)\bar{y}^k + \tau_k\big[\hat{y}^k + \rho_k(K\hat{x}^k + Bw^{k+1} - b)\big]$ by Step~\ref{step:A1-i7}, using the definition \eqref{eq:minimax_form} of $\Lc$, we can easily show that
\myeqn{
\arraycolsep=0.2em
\begin{array}{lcl}
\Lc(x, w, \bar{y}^{k+1}) - (1-\tau_k)\Lc(x, w, \bar{y}^k) &= & \tau_k\iprods{\hat{y}^k + \rho_k(K\hat{x}^k + Bw^{k+1} - b), Kx + Bw - b} \big] \vspace{1ex}\\
&&+ {~} \tau_k F(x, w).
\end{array}
}
Substituting the above two estimates into \eqref{eq:main_bound01}, we can further derive
\myeq{eq:est_1000}{
\arraycolsep=0.2em
\hspace{-2ex}\begin{array}{ll}
\mathbb{E}_{i_k}\big[ \bar{\Lc}_{k+1}(y) & - {~}  \Lc(x, w, \bar{y}^{k+1})  \big] \leq (1-\tau_k)\big[ \bar{\Lc}_k(y)  - \Lc(x, w, \bar{y}^k) \big] \vspace{1.3ex}\\
&+  {~} \frac{1}{2\eta_k}\Exps{i_k}{\norms{\hat{y}^k - y}^2 - \norms{\hat{y}^{k+1} - y}^2 + \norms{\hat{y}^{k+1} - \hat{y}^k}^2} \vspace{1.25ex}\\
&+ {~} \sum_{i=1}^n \frac{\tau_k}{2q_i}\big[\frac{\tau_k\sigma_i}{\tau_0\beta_k}  + (1-q_i)\mu_{f_i} \big] \norms{\tilde{x}_i^k - x_i}^2 \vspace{1ex}\\
&- {~} \Exps{i_k}{ \sum_{i=1}^n \frac{\tau_k}{2q_i}\big(\frac{\tau_k\sigma_i}{\tau_0\beta_k}  + \mu_{f_i} \big)\norms{\tilde{x}_i^{k+1} - x_i}^2 }\vspace{1ex}\\
&- {~}  \frac{\tau_k^2}{2\tau_0^2}\left(\frac{1}{\beta_k} - \rho_k\bar{L}_{\sigma} - L^h_{\sigma} \right)\sum_{i=1}^n\sigma_iq_i\norms{\bar{\tilde{x}}^{k+1}_i - \tilde{x}_i^k}^2\vspace{1.3ex}\\
&- {~} \frac{\rho_k}{2} \norms{K\hat{x}^k + Bw^{k+1} - b - (1-\tau_k)(Kx^k + Bw^k - b)}^2 \vspace{1ex}\\
&- {~}  \frac{(1-\tau_k)}{2}\big[\rho_{k-1} - (1-\tau_k)\rho_k\big] \norms{Kx^k + Bw^{k} - b}^2.
\end{array}
\hspace{-6ex}
}
Next, by \eqref{eq:useful_id3} and $x^{k+1}$ is only different from $\hat{x}^k$ at one block $i_k$, we have
\myeqn{
\arraycolsep=0.2em
\begin{array}{lcl}
\Cc_k &:= & \frac{1}{2\eta_k}\norms{\hat{y}^{k+1} - \hat{y}^k}^2 - \frac{\rho_k}{2} \norms{K\hat{x}^k + Bw^{k+1} - b - (1-\tau_k)(Kx^k + Bw^k - b)}^2 \vspace{1.3ex}\\
& = & \frac{\eta_k}{2}\norms{Kx^{k+1} + Bw^{k+1} - b - (1-\tau_k)(Kx^k + Bw^k - b)}^2 \vspace{1ex}\\
&& - {~} \frac{\rho_k}{2} \norms{K\hat{x}^k + Bw^{k+1} - b - (1-\tau_k)(Kx^k + Bw^{k} - b)}^2 \vspace{0.75ex}\\
&\overset{\tiny\eqref{eq:useful_id3}}{\leq} & \frac{\eta_k\rho_k}{2(\rho_k - \eta_k)}\norms{K(x^{k+1} - \hat{x}^k)}^2  \leq  \frac{\eta_k\rho_k\bar{L}_{\sigma}}{2(\rho_k - \eta_k)}\sum_{i=1}^n\sigma_i\norms{x_i^{k+1} - \hat{x}_i^k}^2.
\end{array}
}
Note also that $\Exps{i_k}{\norms{x_i^{k+1} - \hat{x}^k_i}^2} = \frac{\tau_k^2q_i}{\tau_0^2}\norms{ \bar{\tilde{x}}_i^{k+1} - \tilde{x}^k_i}^2$ due to \eqref{eq:coordinate_update} and Step~\ref{step:A1-i5} of Algorithm~\ref{alg:A1_main}.
Using these expressions, we can estimate
\myeqn{
\arraycolsep=0.2em
\begin{array}{lcl}
\Exps{i_k}{\Cc_k} &\leq & \frac{\eta_k\rho_k\bar{L}_{\sigma}}{2(\rho_k - \eta_k)}\mathbb{E}_{i_k}\Big[ {\displaystyle\sum_{i=1}^n}\sigma_i\norms{x_i^{k+1} - \hat{x}_i^k}^2 \Big]  = \frac{\tau_k^2\rho_k\eta_k}{2\tau_0^2(\rho_k - \eta_k)} \Big[\bar{L}_{\sigma}{\displaystyle\sum_{i=1}^n} q_i\sigma_i\norms{ \bar{\tilde{x}}_i^{k+1} - \tilde{x}^k_i}^2 \Big].
\end{array}
}
Substituting the last inequality into \eqref{eq:est_1000}, we obtain \eqref{eq:est_1001}.
\end{proof}

\beforesubsec
\subsection{Conditions for parameter selection}\label{apdx:le:key_estimate_for_Lyapunov_func}
\aftersubsec
The following lemma provides conditions on the parameters to guarantee a contraction property of the Lyapunov function $\Ec_k(\cdot)$ defined by  \eqref{eq:lyapunov_func}.

\begin{lemma}\label{le:key_estimate_for_Lyapunov_func}
Let $\tau_0$, $\bar{L}_{\sigma}$, and $L^h_{\sigma}$ be defined by \eqref{eq:common_quantities}, and $\set{(x^k, w^k, \bar{y}^k)}$ be generated by  Algorithm~\ref{alg:A1_main}.
Suppose that $\tau_k$, $\beta_k$, $\rho_k$, and $\eta_k$ satisfy the following conditions:
\myeq{eq:param_cond0}{
\arraycolsep=0.3em
\left\{\begin{array}{llcl}
&\rho_{k-1} &\geq & (1-\tau_k)\rho_k, \vspace{1ex}\\
&\eta_k(1-\tau_k) &\geq & \eta_{k-1}, \vspace{1ex}\\
&\frac{\rho_k - \eta_k}{L^h_{\sigma}(\rho_k - \eta_k) + \bar{L}_{\sigma}\rho_k^2}  &\geq & \beta_k, \vspace{1ex}\\
&\frac{\sigma_i\tau_{k-1}^2}{\tau_0\beta_{k-1}} + \mu_{f_i}\tau_{k-1} &\geq &  \frac{1}{(1-\tau_k)} \big[ \frac{\sigma_i\tau_k^2}{\tau_0\beta_k} +  (1-q_i)\mu_{f_i}\tau_k \big], \quad \forall i\in [n].
\end{array}\right.
}
Then, for given $(x, w, y) \in \dom{F}\times\R^d$, the function $\Ec_k(\cdot)$ defined by \eqref{eq:lyapunov_func} satisfies
\myeq{eq:descent_of_lyapunov_func}{
 \Exp{\Ec_{k+1}(x, w, y) }  \leq  (1-\tau_k)\Exp{\Ec_k(x, w, y) }.
}
\end{lemma}

\begin{proof}
From the conditions of \eqref{eq:param_cond0}, we can easily check that
\myeqn{
\frac{1}{\eta_k} \leq \frac{1-\tau_k}{\eta_{k-1}}, \quad  
\rho_{k-1} - (1-\tau_k)\rho_k \geq 0, \quad \text{and}\quad \frac{1}{\beta_k} - \rho_k\bar{L}_{\sigma} - L^h_{\sigma} - \frac{\rho_k\eta_k\bar{L}_{\sigma}}{\rho_k-\eta_k} \geq 0. 
}
Using these relations and the last two conditions of \eqref{eq:param_cond0}, we can simplify \eqref{eq:est_1001} as
\myeq{eq:proof_est20b}{
\hspace{-3ex}\begin{array}{ll}
&\Exps{i_k}{ \bar{\Lc}_{k+1}(y) -  \Lc(x, w, \bar{y}^{k+1}) + \frac{1}{2\eta_k} \norms{\hat{y}^{k+1}  - y }^2  } \vspace{1ex}\\
&\hspace{15ex} \leq (1-\tau_k)\left[ \bar{\Lc}_k(y) - \Lc(x, w, \bar{y}^k)  + \frac{1}{2\eta_{k-1}}\norms{\hat{y}^k - y }^2\right]  \vspace{1ex}\\
&\hspace{15ex} + {~} (1-\tau_k) \sum_{i=1}^n\frac{\tau_{k-1}}{2q_i}\big(\frac{\tau_{k-1}\sigma_i}{\tau_0\beta_{k-1}}  + \mu_{f_i} \big)\norms{\tilde{x}_i^k - x_i }^2 \vspace{1ex}\\
&\hspace{15ex}  - {~}  \Exps{i_k}{\sum_{i=1}^n \frac{\tau_k}{2q_i}\big(\frac{\tau_k\sigma_i}{\tau_0\beta_k}  + \mu_{f_i} \big)\norms{\tilde{x}_i^{k+1} - x_i }^2 }.
\end{array}
\hspace{-3ex}
}
Rearranging this inequality and using $\Ec_k$ defined by \eqref{eq:lyapunov_func},  we obtain
\myeqn{
\begin{array}{ll}
\Exps{i_k}{\Ec_{k+1}(x, w, y)  } &\leq (1-\tau_k)\Ec_k(x, w, y).
\end{array}{\!\!\!\!}
}
Taking full expectation $\Exp{\cdot}$ given $(x, w, y)$ on both sides of the last inequality so that $\Exp{\Exps{i_k}{\cdot}  } = \Exp{\Exp{\cdot \mid \Fc_{k}} } = \Exp{\cdot  }$, we eventually get
\myeqn{
\Exp{\Ec_{k+1}(x, w, y)  } \leq (1-\tau_k)\Exp{\Ec_k(x, w, y) },
}
which proves \eqref{eq:descent_of_lyapunov_func}.
\end{proof}

\beforesubsec
\subsection{The proof of Theorem~\ref{th:convergence_rate}: General convex case}\label{apdx:le:convergence_result}
\aftersubsec
Since $\mu_{f_i} = 0$ for all $i \in [n]$, if we assume that the conditions of \eqref{eq:param_cond0} are tight, 
then we can easily derive that
\myeq{eq:update_param}{
\rho_k := \frac{\rho_{k-1}}{1-\tau_k}  \qquad  \text{and} \qquad \tau_k := \frac{\tau_{k-1}}{\tau_{k-1} + 1},
}
where $\rho_0 > 0$ is given and $\tau_0$ is defined by \eqref{eq:common_quantities}.
Let us also update $\eta_k$ as $\eta_k := \frac{\rho_k}{2}$.
Then, it is straightforward to prove that
\myeq{eq:tau_update}{
\arraycolsep=0.2em
\begin{array}{l}
\tau_k := \frac{\tau_0}{\tau_0k + 1}, \ \ \rho_k := \rho_0(\tau_0 k + 1), \ \ \eta_k := \frac{\rho_0(\tau_0 k + 1)}{2}, \ \  \omega_k := \prod_{i=1}^k(1-\tau_i) = \frac{1}{\tau_0k + 1}.
\end{array}
}
Moreover, the third condition $\frac{\rho_k - \eta_k}{L^h_{\sigma}(\rho_k - \eta_k) + \bar{L}_{\sigma}\rho_k^2}  \geq \beta_k$ of \eqref{eq:param_cond0}  becomes $\frac{1}{L^h_{\sigma} + 2\bar{L}_{\sigma}\rho_k} \geq \beta_k$.
Hence, we can update $\beta_k$ as following to guarantee this condition: 
\myeq{eq:betagamma_update}{
\beta_k := \frac{1}{L^h_{\sigma} + 2\bar{L}_{\sigma}\rho_k} = \frac{1}{L^h_{\sigma} + 2\bar{L}_{\sigma}\rho_0(\tau_0 k + 1)}.
}
In summary, it is clear that the update rule \eqref{eq:param_update1} satisfies all the conditions of \eqref{eq:param_cond0}.

Next, from \eqref{eq:descent_of_lyapunov_func} and \eqref{eq:tau_update}, we can show that
\myeq{eq:lm3_proof1}{
\arraycolsep=0.3em
\begin{array}{lcl}
\Exp{\Ec_{k+1}(x, w, y)  } & \overset{\tiny\eqref{eq:descent_of_lyapunov_func}}{\leq} & \Big[ \prod_{i=1}^k(1-\tau_i) \Big]\Exp{\Ec_1(x, w, y)  }  \vspace{1ex} \\
&  \overset{\tiny\eqref{eq:tau_update}}{=}  & \frac{1}{\tau_0k + 1}\Exp{\Ec_1(x, w, y) }.
\end{array}
}
Using \eqref{eq:est_1001} and the definition \eqref{eq:lyapunov_func} of $\Ec_{k}$, we have
\myeq{eq:E1_est}{
\hspace{-3ex}
\arraycolsep=0.1em
\begin{array}{lcl}
\Ec_1(x, w, y) &\overset{\eqref{eq:lyapunov_func}}{=} & \bar{\Lc}_1(y) - \Lc(x,w,\bar{y}_1) + \frac{1}{2\eta_0}\norms{y-\hat{y}^1}^2 + \sum_{i=1}^n\frac{\tau_{0}\sigma_i}{2\beta_0q_i}\norms{\tilde{x}_i^{1} - x_i}^2 \vspace{1ex}\\
&\overset{\eqref{eq:est_1001}}{\leq} & (1-\tau_0)[\bar{\Lc}_0(y) - \Lc(x,w,\bar{y}_0)] + \sum_{i=1}^n\frac{\tau_0\sigma_i}{2\beta_0q_i} \norms{\tilde{x}_i^0 - x_i}^2 +  \frac{1}{2\eta_0}\norms{y - \hat{y}^0}^2 \vspace{1ex}\\
& \leq & \hat{\Ec}_0(x, w, y), \quad (\text{since $1-\tau_0\leq 1$})
\end{array}
\hspace{-2ex}
}
where $\hat{\Ec}_0(x, w, y)$ is defined as
\myeqn{
\arraycolsep=0.2em
\begin{array}{lcl}
\hat{\Ec}_0(x, w, y) &:=&  F(x^0, w^0) - \Lc(x, w, \hat{y}^0)  + \iprods{y, Kx^0 + Bw^0 - b}  + \frac{1}{2\eta_0}\norms{\hat{y}^0 - y}^2 \vspace{1.5ex}\\
&& + {~}  \frac{\rho_{-1}}{2}\norms{Kx^0 + Bw^0 - b}^2 + \frac{(L^h_{\sigma} + 2\rho_0\bar{L}_{\sigma})\tau_0}{2}\norms{x - x^0}_{\sigma/q}^2.
\end{array}
}
Denoting $u^0 := Kx^0 + Bw^0 - b$.
Then, since $\rho_{-1} = \rho_0(1-\tau_0)$ and $2\eta_0 = \rho_0$, we have
\rvtext{
\myeqn{
\arraycolsep=0.2em
\begin{array}{lcl}
\Tc_{[3]} &:= & \iprods{y, Kx^0 + Bw^0 - b}  + \frac{1}{2\eta_0}\norms{\hat{y}^0 - y}^2 + \frac{\rho_{-1}}{2}\norms{Kx^0 + Bw^0 - b}^2 \vspace{1ex}\\
&\leq & \frac{1}{2\rho_0}\norms{y}^2 + \frac{\rho_0}{2}\norms{u^0}^2 + \frac{1}{\rho_0}\norms{\hat{y}^0 - y}^2 + \frac{\rho_0(1-\tau_0)}{2}\norms{u^0}^2 \vspace{1ex}\\
&\leq & \frac{1}{\rho_0}\norms{\hat{y}^0}^2 + \frac{2}{\rho_0}\norms{\hat{y}^0 - y}^2 + \frac{\rho_0(2-\tau_0)}{2}\norms{u^0}^2,
\end{array}
}
where the last inequality comes from $\frac{1}{2}\norms{y}^2 \leq \norms{\hat{y}^0}^2 + \norms{y - \hat{y}^0}^2$.
}
Using $\Tc_{[3]}$, and the fact that $\hat{y}^0 := 0$ and  $-\Lc(x, w, \hat{y}^0) \leq -D(\hat{y}^0)$, we can further simplify $\hat{\Ec}_0$ as
\rvtext{
\myeq{eq:E0_hat}{
\arraycolsep=0.2em
\begin{array}{lcl}
\hat{\Ec}_0(x, w, y) & \leq & \hat{\Ec}_0(x, y) :=  F(x^0, w^0) - D(\hat{y}^0)  + \frac{2}{\rho_0}\norms{ y - \hat{y}^0}^2 + \frac{1}{\rho_0}\norms{\hat{y}^0}^2 \vspace{1.5ex}\\
&& + {~}  \frac{\rho_0(2-\tau_0)}{2}\norms{Kx^0 + Bw^0 - b}^2 + \frac{(L^h_{\sigma} + 2\rho_0\bar{L}_{\sigma})\tau_0}{2}\norms{x - x^0}_{\sigma/q}^2.
\end{array}
}
}
Now, by convexity of $f$, using \eqref{eq:non_smooth} and \eqref{eq:x_r_breve} we can show that
\myeqn{
f(x^k) \overset{\tiny\eqref{eq:x_r_breve}}{=}  f\left(\sum_{l=0}^k\gamma_{k,l}\tilde{x}^l\right) \leq \sum_{l=0}^k\gamma_{k,l}f(\tilde{x}^l)  \overset{\tiny\eqref{eq:non_smooth}}{=}  \bar{f}^k.
}
Therefore, we can derive
\myeqn{
\arraycolsep=0.1em
\begin{array}{ll}
\Lc_{\rvtext{\rho_{k-1}}}(x^{k}, w^{k}, y) - \Lc(x, w, \bar{y}^{k}) & =  F(x^{k}, w^{k})   + \psi_{\rho_{k-1}}(x^{k}, w^{k}, y)  - \Lc(x, w, \bar{y}^{k}) \vspace{1ex}\\
&\leq   \bar{f}^{k} + h(x^{k}) + g(w^{k})  + \psi_{\rho_k}(x^{k},  w^{k}, y) - \Lc(x, w, \bar{y}^{k})  \vspace{1ex}\\
&\overset{\tiny\eqref{eq:lyapunov_func}}{\leq}   \Ec_{k}(x, w, y).
\end{array}
}
Combining this inequality, $\mathbb{E}[\hat{\Ec}_0(x, y)] = \hat{\Ec}_0(x, y)$, \eqref{eq:lm3_proof1},  \eqref{eq:E1_est}, and \eqref{eq:E0_hat}, we get
\myeq{eq:key_estimate3}{
\Exp{\Lc_{\rho_{k-1}}(x^{k}, w^k, y) - \Lc(x, w, \bar{y}^{k})  }  \leq \frac{\hat{\Ec}_0(x, y)}{\tau_0k + 1 - \tau_0}.
}
Since $\Lc(x^k, w^k, y) \leq \Lc_{\rho_{k-1}}(x^k, w^k, y)$, \eqref{eq:key_estimate3} leads to 
\myeqn{
\arraycolsep=0.1em
\begin{array}{lcl}
\Gc_{\Zc}(x^k, w^k, \bar{y}^k) &\overset{\tiny\eqref{eq:gap_func}}{=} & \displaystyle\max_{(x,w,y)\in\Zc}\Exp{\Lc(x^k, w^k, y) - \Lc(x, w, \bar{y}^k)} \vspace{1ex}\\
&  \leq &  \frac{1}{\tau_0k + 1 - \tau_0}\Big[F(x^0, w^0) - D(\hat{y}^0) + 
\rvtext{\frac{\rho_0(2-\tau_0)}{2}\norms{Kx^0 + Bw^0 - b}^2 }\vspace{1ex}\\
&& + {~} \rvtext{ \frac{1}{\rho_0}\norms{\hat{y}^0}^2 } + {\displaystyle\sup_{(x,y)\in\Xc\times\Yc}} \left\{ \rvtext{ \tfrac{2}{\rho_0}\norms{ y - \hat{y}^0}^2 } +  \tfrac{(L^h_{\sigma} + 2\rho_0\bar{L}_{\sigma})\tau_0}{2}\norms{x - x^0}_{\sigma/q}^2 \right\}\Big],
\end{array}
}
which proves the last inequality of \eqref{eq:residual_bound}.

Next, using the saddle-point condition \eqref{eq:saddle_point}, we can show that 
\myeqn{
F^{\star} = F(x^{\star}, w^{\star}) = \Lc(x^{\star}, w^{\star}, \bar{y}^k) \overset{\tiny\eqref{eq:saddle_point}}{\leq} \Lc(x^k, w^k, y^{\star}) = F(x^k, w^k) + \iprods{y^{\star}, Kx^k + Bw^k - b}.
}
This implies that $\Exp{F(x^k, w^k) - F^{\star} + \iprods{y^{\star}, Kx^k + Bw^k - b}} \geq 0$.
Here $\Exp{\cdot} = \Exp{\cdot }$.
On the other hand, using \eqref{eq:key_estimate3}, we also have
\myeq{eq:lm3_proof3}{
\hspace{-2ex}
\Exp{ F(x^k, w^k) \! - \! F^{\star} \! + \! \iprods{y^{\star}, Kx^k \! + \! Bw^k \! - \! b} + \frac{\rho_{k-1}}{2}\norms{Kx^k \! + \! Bw^k \! - \! b}^2 } \leq \frac{\hat{\Ec}_0(x^{\star}, y^{\star})}{\tau_0k \! + \! 1 \! -\! \tau_0}.
\hspace{-2ex}
}
Hence, together with $\rho_{k-1} = \rho_0(\tau_0k + 1 - \tau_0)$, we obtain 
\myeqn{
\Exp{\norms{Kx^k + Bw^k - b}^2} \leq  \frac{2\hat{\Ec}_0(x^{\star}, y^{\star})}{\rho_0(\tau_0k + 1 - \tau_0)^2}.
}
Moreover, from \eqref{eq:Ebar_0}, we have $\bar{\Ec}_0 = \hat{\Ec}_0(x^{\star},  y^{\star})$.
Thus \eqref{eq:lm3_proof3} implies
\myeqn{
\begin{array}{lcl}
\left\vert \Exp{ F(x^k, w^k) - F^{\star} } \right\vert &\leq & \frac{1}{\tau_0k + 1 - \tau_0}\bar{\Ec}_0 + \norms{y^{\star}}\left(\Exp{\norms{Kx^k + Bw^k - b}^2}\right)^{1/2} \vspace{1ex}\\
&\leq & \frac{1}{\tau_0k + 1 - \tau_0}\left[\bar{\Ec}_0   + \norms{y^{\star}}\left(\frac{2}{\rho_0}\bar{\Ec}_0\right)^{1/2}\right],
\end{array}
}
which proves the first two lines of \eqref{eq:residual_bound}.

\rvtext{Now, let $D(\bar{y}^k) := \min_{x, w}\Lc(x, w, \bar{y}^k)$ be the dual function.
Then, we have
\myeqn{
\arraycolsep=0.3em
\begin{array}{lcl}
D(\bar{y}^k) &= & {\displaystyle\min_{x, w}} \big\{\phi(x) + g(w) + \iprod{Kx + Bw - b, \bar{y}^k} \big\} \vspace{1ex}\\
& = & -\phi^{\ast}(-K^{\top}\bar{y}^k) - g^{\ast}(-B^{\top}\bar{y}^k) - \iprods{b, \bar{y}^k}.
\end{array}
}
Therefore,  $\dom{D} = \sets{ y \in \R^d \mid  -K^{\top}y \in \dom{\phi^{\ast}}, \  -B^{\top}y \in \dom{g^{\ast}}}$. 
Let us show that $\bar{y}^k \in \dom{D}$.
Firstly, by the assumption that $\phi^{\ast}$ is $M_{\phi^{\ast}}$-Lipschitz continuous, we have  $\dom{\phi^{\ast}} = \R^d$. 
Hence, we only need to prove $-B^{\top}\bar{y}^k \in \dom{g^{\ast}}$. 
Indeed, from \eqref{eq:w_step1}, we have $0 \in \partial g(w^{k+1}) + B^{\top}(\hat{y}^k + \rho_k(Bw^{k+1} + K\hat{x}^k - b))$,  which becomes
\myeqn{
w^{k+1} \in \partial g^{\ast}(-B^{\top}(\hat{y}^k + \rho_k(Bw^{k+1} + K\hat{x}^k - b))).
}
Thus $-B^{\top}(\hat{y}^k + \rho_k(Bw^{k+1} + K\hat{x}^k - b)) \in \dom{g^{\ast}}$. 
Also from Sep~\ref{step:A1-i7} of our Algorithm~\ref{alg:A1_main}, $\bar{y}^{k+1} = (1 - \tau_k)\bar{y}^k + \tau_k[\hat{y}^k + \rho_k(Bw^{k+1} + K\hat{x}^k - b)]$, where $\tau_k \in (0,1)$. 
Therefore, if $-B^{\top}\bar{y}^k \in \dom{g^{\ast}}$, then $-B^{\top}\bar{y}^{k+1} \in \dom{g^{\ast}}$. 
As a result, if we assume that $-B^{\top}\bar{y}^{0} \in \dom{g^{\ast}}$ at the initialization of Algorithm~\ref{alg:A1_main}, then  $-B^{\top}\bar{y}^{k} \in \dom{g^{\ast}}$ for all $k\geq 0$.
We conclude that $\bar{y}^k \in \dom{D}$, which implies $D(\bar{y}^k) > -\infty$ for all $k\geq 0$.
Finally, since $\phi$ and $g$ are proper, by the Fenchel-Moreau theorem \cite[Theorem 13.37]{Bauschke2011}, $\phi^{*}$ and $g^{*}$ are also proper, and hence $D$ is proper. 
This shows that $D(y) < +\infty$ for all $y \in \dom{D}$.
Combining two cases, we can state that $D(\bar{y}^k)$ is finite.
}

Next, we note that $\bar{y}^k$ is independent of $i_{k-1}$, and therefore, independent of $(x^k, w^k)$.
Hence, from $D(\bar{y}^k) = \min_{x, w}\Lc(x, w, \bar{y}^k)$, using the optimality condition of this minimization problem, for any $\bar{w}^{k-1} \in \partial{g}^{*}(-B^{\top}\bar{y}^k)$ and any $\bar{x}^{k-1} \in \partial{\phi}^{*}(-K^{\top}\bar{y}^k)$, we have $D(\bar{y}^k) = \Lc(\bar{x}^k, \bar{w}^{k-1}, \bar{y}^{k-1})$, where $\phi^{*}$ is the Fenchel conjugate of $\phi := f + h$, and we shift the index to $k-1$ to show that both $\bar{x}^{k-1}$ and $\bar{w}^{k-1}$ are $\Fc_{k-1}$-measurable, and independent of $(x^k, w^k)$.
Therefore, one has
\myeqn{
\arraycolsep=0.3em
\begin{array}{lcl}
D^{\star} - D(\bar{y}^k) & \overset{\tiny\eqref{eq:saddle_point}}{\leq} & \Lc(x^k, w^k, y^{\star}) - \Lc(\bar{x}^{k-1}, \bar{w}^{k-1}, \bar{y}^k) \vspace{1ex}\\ 
& \leq & \Lc_{\rho_{k-1}}(x^k, w^k, y^{\star}) - \Lc(\bar{x}^{k-1}, \bar{w}^{k-1}, \bar{y}^k).
\end{array}
}
Here, we have used $D^{\star} = F^{\star} \overset{\eqref{eq:saddle_point}}{\leq} \Lc(x^k, w^k, y^{\star}) \leq \Lc_{\rho_{k-1}}(x^k, w^k, y^{\star})$.
From \eqref{eq:descent_of_lyapunov_func}, by induction and similar to the proof of \eqref{eq:key_estimate3}, we get
\myeqn{
\mathbb{E}\big[ \Lc_{\rho_{k-1}}(x^{k}, w^k, y^{\star}) - \Lc(x, w, \bar{y}^{k}) \big] \leq \mathbb{E}\big[\Ec_k(x, w, y^{\star}) \big]  \leq \frac{\hat{\Ec}_0(x, y^{\star})}{\tau_0k + 1 - \tau_0},
}
where $\hat{\Ec}_0$ is defined by \eqref{eq:E0_hat}.
Now, we substitute $x = \bar{x}^{k-1}$ and $w = \bar{w}^{k-1}$ into the last inequality to get
\myeqn{
\mathbb{E}\big[ \Lc_{\rho_{k-1}}(x^{k}, w^k, y^{\star}) - \Lc(\bar{x}^{k-1}, \bar{w}^{k-1}, \bar{y}^{k}) \big]  \leq \frac{\hat{\Ec}_0(\bar{x}^{k-1}, y^{\star})}{\tau_0k + 1 - \tau_0}.
}
Here, the expectation is now conditioned on  $(\bar{x}^{k-1}, \bar{w}^{k-1})$, which is random.
Taking the full expectation, and combine the result with the above estimate, we get
\myeqn{
\mathbb{E}\big[ D^{\star} - D(\bar{y}^k) \big] \leq \frac{\mathbb{E}[\hat{\Ec}_0(\bar{x}^{k-1}, y^{\star})]}{\tau_0k + 1 - \tau_0}.
}
In addition, since $\phi^{*}$ is $M_{\phi^{*}}$-Lipschitz continuous, almost surely, we have
\myeqn{
\sup\big\{ \norms{\bar{x}^{k-1} - x^0}^2_{\sigma/q} : \bar{x}^{k-1} \in \partial{\phi}^{*}(-K^{\top}\bar{y}^k) \big\} \leq M_0 := \sup_{ \norms{\bar{x}^{k-1}} \leq M_{\phi^{*}} } \big\{ \norms{\bar{x}^{k-1} - x^0}^2_{\sigma/q} \big\}. 
}
The last two inequalities lead to
\myeqn{
\arraycolsep=0.1em
\begin{array}{lcl}
\Exp{D^{\star} - D(\bar{y}^k)}
&\overset{\tiny\eqref{eq:E0_hat}}{\leq} & \frac{1}{\tau_0k + 1 - \tau_0}\Big[ F(x^0,w^0) - D(\hat{y}^0) + \frac{2}{\rho_0}\norms{\hat{y}^0 - y^{\star}}^2 + \frac{1}{\rho_0}\norms{\hat{y}^0}^2  \vspace{0.5ex}\\
&& + {~} \frac{\rho_0(2 - \tau_0)}{2}\norms{Kx^0 + Bw^0 - b}^2 + \frac{(L^h_{\sigma} + 2\rho_0\bar{L}_{\sigma})\tau_0M_0}{2} \Big]  =  \frac{\bar{\Fc}_0}{\tau_0k + 1 - \tau_0},
\end{array}
\vspace{-2ex}
}
which proves the third line of \eqref{eq:residual_bound}.
\Eproof

\beforesubsec
\subsection{The proof of Theorem~\ref{th:convergence_rate2}: Strongly convex case}\label{apdx:th:convergence_rate2}
\aftersubsec
We first show that if $\tau_k$, $\rho_k$, $\eta_k$, and $\beta_k$ are updated as in Theorem~\ref{th:convergence_rate2}, then they  satisfy the conditions of \eqref{eq:param_cond0}.
First, it is obvious to show that $\rho_k$, $\eta_k$, and $\beta_k$ satisfy the first three conditions of \eqref{eq:param_cond0}.
Next, since $\tau_k$  is updated as in Theorem~\ref{th:convergence_rate2}, it satisfies $1-\tau_k = \frac{\tau_k^2}{\tau_{k-1}^2}$.
Hence, we obtain 
\myeqn{
\frac{\tau_0}{\tau_0k+1} \leq \tau_k \leq \frac{2\tau_0}{\tau_0k+2}, ~~~~~~\prod_{i=1}^k(1-\tau_i) = \frac{\tau_k^2}{\tau_0^2} \leq \frac{4}{(\tau_0k+2)^2}, ~~~~\text{and}~~~~\rho_k = \frac{\tau_{k-1}^2}{\tau_k^2}\rho_{k-1}.
}
Then, by induction, we get $\rho_k = \frac{\rho_0\tau_0^2}{\tau_k^2}$.
Therefore, $\beta_k = \frac{\tau_k^2}{L^h_{\sigma}\tau_k^2 + 2\bar{L}_{\sigma}\rho_0\tau_0^2}$.
Consequently, one can show that $\frac{\tau_k}{\tau_0\beta_k} - \frac{\tau_k}{\tau_0\beta_{k-1}} = 2 \bar{L}_{\sigma}\rho_0\tau_0$.

Now, we verify the last condition of \eqref{eq:param_cond0}.
Multiplying this condition by $\frac{1-\tau_k}{\tau_k} = \frac{\tau_k}{\tau_{k-1}^2}$, it is equivalent to
\myeqn{
\frac{\mu_{f_i}}{\sigma_i}\left[ \frac{\tau_{k}}{\tau_{k-1}} - (1- q_i)\right] \geq \frac{\tau_k}{\tau_0\beta_k} - \frac{\tau_k}{\tau_0\beta_{k-1}} = 2\bar{L}_{\sigma}\rho_0\tau_0.
}
It is easy to check that $\frac{\tau_{k}}{\tau_{k-1}} = \sqrt{1 - \tau_k}$ is increasing. 
Using $\tau_0 \leq q_i$ from \eqref{eq:common_quantities}, the above inequality holds if
\myeqn{
\begin{array}{ll}
& 2\bar{L}_{\sigma}\rho_0\tau_0 \leq \frac{\mu_{f_i}}{\sigma_i}\left[ \frac{\tau_{1}}{\tau_{0}} - (1- \tau_0) \right] = \frac{\mu_{f_i}}{\sigma_i}\left[ \frac{\tau_{1}}{\tau_{0}} + \tau_0 - 1\right],
\end{array}
}
which is equivalent to $0 < \rho_0 \leq \min_{i \in[n]}\set{ \frac{\mu_{f_i}}{\sigma_i\bar{L}_{\sigma}}}\frac{\sqrt{\tau_0^2 + 4} + \tau_0 - 2}{4\tau_0}$.
Using $\frac{\sqrt{\tau_0^2 + 4} + \tau_0 - 2}{4\tau_0} \geq \frac{1}{4}$, we can simplify this expression by a tighter one $0 < \rho_0 \leq \frac{1}{4\bar{L}_{\sigma}}\min_{i\in [n]}\frac{\mu_{f_i}}{\sigma_i}$, which is exactly the last condition in Theorem~\ref{th:convergence_rate2}.

The remaining proof of Theorem~\ref{th:convergence_rate2} is similar to the proof of Theorem~\ref{th:convergence_rate} but using $\tilde{\Ec}_0(x, y)$ instead of $\hat{\Ec}_0(x, y)$ in the key bound \eqref{eq:key_estimate3}, where
\myeqn{
\arraycolsep=0.2em
\begin{array}{lcl}
\tilde{\Ec}_0(x, y) &:=&  F(x^0, w^0) - D(\hat{y}^0)  + 
\rvtext{\frac{2}{\rho_0}\norms{ y - \hat{y}^0}^2 + \frac{1}{\rho_0}\norms{\hat{y}^0}^2 + \frac{\rho_0(2-\tau_0)}{2}\norms{Kx^0 + Bw^0 - b}^2} \vspace{1.5ex}\\
&& + {~} \frac{\tau_0}{2}\sum_{i=1}^n\frac{1}{q_i}\big[(L^h_{\sigma} + 2\rho_0\bar{L}_{\sigma})\sigma_i  + \mu_{f_i} \big] \norms{x_i - x^0_i}^2.
\end{array}
}
\rvtext{Similarly, we also replace $\bar{\Fc}_0$ and $\bar{R}_\Zc^2$ by $\tilde{\Fc}_0$ and $\tilde{R}_\Zc^2$, respectively}.
To avoid repetition, we omit the detailed derivation here.
\Eproof

\beforesec
\section{Conclusions}\label{sec:conclusion}
\aftersec
We have developed a unified randomized block-coordinate alternating primal-dual algorithm to solve a generic class of nonsmooth and constrained convex optimization problems of the form \eqref{eq:primal_form} and its dual problem \eqref{eq:dual_form}.
Our algorithm is new and achieves optimal convergence rates (up to a constant factor) for both general convex case and strongly convex case.
Our rates are on three criteria and non-ergodic on the primal sequence.
We have also specified our algorithm to handle two special cases, which commonly appear in the literature.
This leads to new variants where our convergence rates guarantees are still applied.
We have tested our algorithm on two well-studied examples and compared it with two state-of-the-art algorithms.
We have observed that our algorithm has encouraging performance on different experiments of real and synthetic datasets.
Our next step is to extend this approach to convex optimization models with nonlinear constraints and general convex-concave saddle-point problems.

\vskip 2mm
\noindent{\bf Acknowledgments.}
This work is partly supported by the Office of Naval Research under Grant No.  ONR-N00014-20-1-2088 (2020 - 2023), and the Nafosted Vietnam, Grant No. 101.01-2020.06 (2020 - 2022).

\bibliographystyle{plain}

\end{document}